\documentclass[10pt]{article}
\usepackage{rotating}
\usepackage{amsmath,amssymb}
\usepackage{graphics}
\usepackage[all]{xy}
\xyoption{v2}
\xyoption{knot}

\newtheorem{theorem}{Theorem}
\newtheorem{definition}{Definition}

\newtheorem{example}{Example}
\newcommand{\MySection}[1]
{\section{ #1}}

\newcommand{\bcal}[1]{\mbox{\boldmath${\cal {#1}}$}}

\newenvironment{fake}{\relax}{\relax} 

\newenvironment{proof}
{\noindent{\bf Proof}}

\begin{document}


\title{Enrichment as Categorical Delooping I: Enrichment Over Iterated Monoidal Categories}
\author{Stefan Forcey}
\maketitle
\begin{abstract}
Joyal and Street note in their paper on braided monoidal categories \cite{JS}
that the 2--category ${\cal V}$--Cat of categories enriched over a braided monoidal category 
   ${\cal V}$ is not itself braided in 
   any way that is based upon the braiding of ${\cal V}$. What is meant by ``based upon'' here will be made more 
   clear in the present paper.
   The exception that they mention is the case in which ${\cal V}$ is symmetric, which leads to 
   ${\cal V}$--Cat being symmetric as well.
   The symmetry in ${\cal V}$--Cat is based upon the symmetry of ${\cal V}$. 
   The motivation behind this paper is in part to 
   describe
   how these facts relating ${\cal V}$ and ${\cal V}$--Cat are in turn related to a categorical analogue of topological 
  delooping first mentioned 
   by Baez and Dolan in \cite{Baez1}.
   To do so 
   I need to pass to a more general setting than braided and symmetric categories -- in 
   fact the $k$--fold monoidal categories 
   of Balteanu et al in \cite{Balt}. It seems that the analogy of loop spaces is a good guide for how to define the concept of enrichment
   over various types of monoidal objects, including $k$--fold monoidal 
   categories and their higher dimensional counterparts.
   The main result is that for ${\cal V}$ a $k$--fold monoidal category, ${\cal V}$--Cat 
   becomes a $(k-1)$--fold monoidal $2$--category 
   in a canonical way. I indicate how this process may be iterated by enriching over ${\cal V}$--Cat, along the way defining the 
   3--category of categories enriched over ${\cal V}$--Cat.
   In the next paper I hope to make precise the $n$--dimensional case and
   to show how the group completion of the nerve of ${\cal V}$ is related to the 
loop space of the group completion of the nerve of ${\cal V}$--Cat.
\end{abstract}
\begin{fake}\end{fake}

\SloppyCurves{
      
      \clearpage
            \MySection{Introduction}



    A major goal of higher dimensional category theory is to discover ways of exploiting the 
    connections between homotopy coherence and categorical coherence. Stasheff \cite{Sta} and
    MacLane \cite{Mac} showed that monoidal categories are precisely analogous to 1--fold
    loop spaces. There is a similar connection between symmetric monoidal categories and
    infinite loop spaces. The first step in filling in the gap between 1 and infinity was made in
    \cite{ZF} where it is shown that the group completion of the nerve of a braided monoidal category
    is a 2--fold loop space.   
    In \cite{Balt} the authors finished this process by, in their words, ``pursuing an analogy to the tautology 
    that an $n$--fold loop space is a loop
    space in the category of $(n-1)$--fold loop spaces.'' The first
    thing they focus on is the fact that a braided category is a special case of a carefully
    defined 2--fold monoidal category. Based on their observation of the  correspondence between
    loop spaces and monoidal categories, they iteratively define the notion of $n$--fold
    monoidal category as a monoid in the category of $(n-1)$--fold monoidal categories.
    In their view
    ``monoidal'' functors should be defined in a more ``lax'' way than is usual in order to avoid
    strict commutativity of 2--fold and higher monoidal categories. In \cite{Balt} a symmetric category
    is seen as a category that is $n$--fold monoidal for all $n$.
    
    The main result in \cite{Balt} is that their definition of iterated monoidal categories
    exactly corresponds to $n$--fold loop spaces for all $n$. They show that the group
    completion of the nerve of an $n$--fold monoidal category is an $n$--fold loop space. Then they describe 
    an operad in the category of small categories which parameterizes the
    algebraic structure of an iterated monoidal category.  They show that the nerve of
    this categorical operad is a topological operad which is equivalent
    to the little $n$--cubes operad. This latter operad, as shown in \cite{BV1} and \cite{May},
    characterizes the notion of $n$--fold loop space.  Thus the main result in \cite{Balt} is a categorical
    characterization of $n$--fold loop spaces.  
    
    The present paper pursues the hints of a categorical delooping that are suggested by the facts that
    for a symmetric category, the 2--category of categories enriched over it is again symmetric, while for a 
    braided category the 2--category of categories enriched over it is merely monoidal. Section 1 reviews enrichment
    and Section 2 investigates just what obstacles arise when defining a product based on a braiding and 
    attempting to define further a braiding of that derived product. Section 3 goes over the recursive definition
    of the $k$--fold monoidal categories of \cite{Balt}, altered here to include a coherent associator. 
    The immediate question is whether the delooping phenomenon happens in general for these $k$--fold 
    monoidal categories.  The answer is yes, once enriching over a $k$--fold monoidal category is 
    carefully defined in Section 4. The definition
    also provides for iterated delooping, and all the information included in the axioms for the $k$--fold category
    is exhausted in the process, as described in Section 5. 
    It seems that passing to the category of enriched categories basically reduces the number
    of products so that for ${\cal V}$ a $k$--fold monoidal $n$--category, ${\cal V}$--Cat becomes a $(k-1)$--fold 
    monoidal $(n+1)$--category. This picture was anticipated by Baez and Dolan \cite{Baez1} in the context where the 
    $k$--fold monoidal $n$--category is specifically a (weak) $(n+k)$--category with only one object,
    one 1--cell, etc. up to only one $k$--cell. The construction of a $k$--fold monoidal $n$--category
    here is a bit different, from the bottom up as it were. The next question is whether and how the two constructions
    overlap. Since the enriched category construction gives strict $n$--categories, it seems to be a special case of
    what they expect to be true in general. 
    
    In \cite{StAlg} Street defines the nerve of a strict n-category.
    Recently Duskin in \cite{Dusk} has worked out the description of the nerve of a bicategory. 
    This allows us to ask whether
    these nerves will prove to be the logical link to loop spaces for higher dimensional iterated monoidal categories. 
    The most basic statement
    should be that the loop space of the group completion of the nerve of ${\cal V}$--Cat is precisely the group completion of 
    the nerve of ${\cal V}$. I
    will give examples that indicate that this is not necessarily the case, and hopefully more elucidation of
    the question as well, in further work.

    I have organized the paper so that sections can largely stand alone, so please skip them when able, 
    and forgive redundancy when it occurs.


    \clearpage
    \newpage
    \MySection{Review of Categories Enriched Over a Monoidal Category}
    In this section I briefly review the definition of a category enriched over a monoidal category ${\cal V}$. 
  Enriched functors and 
    enriched natural transformations make the collection of enriched categories into a 2-category ${\cal V}$-Cat. 
    This section is not meant to be complete. It is included due to, and its contents determined by, how often the 
    definitions herein are referred to and followed as models in the rest of the paper. The definitions and proofs
    can be found in more or less detail in \cite{Kelly} and \cite{EK1} and of course
    in \cite{MacLane}.
    
    \begin{definition} For our purposes a {\it monoidal category} is a category ${\cal V}$
      together with a functor 
      $\otimes: {\cal V}\times{\cal V}\to{\cal V}$  and an object $I$ such that
      \begin{enumerate}
      \item $\otimes$ is  associative up to the coherent natural transformations $\alpha$. The coherence
      axiom is given by the commuting pentagon
      
      \noindent
      	\begin{center}
      	\resizebox{5.5in}{!}{
              $$
              \xymatrix@C=-25pt{
              &((U\otimes V)\otimes W)\otimes X \text{ }\text{ }
              \ar[rr]^{ \alpha_{UVW}\otimes 1_{X}}
              \ar[ddl]^{ \alpha_{(U\otimes V)WX}}
              &&\text{ }\text{ }(U\otimes (V\otimes W))\otimes X
              \ar[ddr]^{ \alpha_{U(V\otimes W)X}}&\\\\
              (U\otimes V)\otimes (W\otimes X)
              \ar[ddrr]|{ \alpha_{UV(W\otimes X)}}
              &&&&U\otimes ((V\otimes W)\otimes X)
              \ar[ddll]|{ 1_{U}\otimes \alpha_{VWX}}
              \\\\&&U\otimes (V\otimes (W\otimes X))&&&
             }
             $$
             }
	                    \end{center}
	                    
    \item $I$ is a strict $2$-sided unit for $\otimes$.
    \end{enumerate}
    \end{definition}
    \begin {definition} \label{V:Cat} A (small) ${\cal V}$ {\it--Category} ${\cal A}$ is a set $\left|{\cal A}\right|$ of 
    {\it objects}, 
    a {\it hom-object} ${\cal A}(A,B) \in \left|{\cal V}\right|$ for
    each pair of objects of ${\cal A}$, a family of {\it composition morphisms} $M_{ABC}:{\cal A}(B,C)
    \otimes{\cal A}(A,B)\to{\cal A}(A,C)$ for each triple of objects, and an {\it identity element} $j_{A}:I\to{\cal A}(A,A)$ for each object.
    The composition morphisms are subject to the associativity axiom which states that the following pentagon commutes

          \noindent
          	          \begin{center}
    	          \resizebox{5.5in}{!}{
          $$
          \xymatrix@C=-15pt{
          &({\cal A}(C,D)\otimes {\cal A}(B,C))\otimes {\cal A}(A,B)\text{ }\text{ }
          \ar[rr]^{\scriptstyle \alpha}
          \ar[dl]^{\scriptstyle M \otimes 1}
          &&\text{ }\text{ }{\cal A}(C,D)\otimes ({\cal A}(B,C)\otimes {\cal A}(A,B))
          \ar[dr]^{\scriptstyle 1 \otimes M}&\\
          {\cal A}(B,D)\otimes {\cal A}(A,B)
          \ar[drr]^{\scriptstyle M}
          &&&&{\cal A}(C,D)\otimes {\cal A}(A,C)
          \ar[dll]^{\scriptstyle M}
          \\&&{\cal A}(A,D))&&&
          }$$
          }
         	                    \end{center}
     	                    
    and to the unit axioms which state that both the triangles in the following diagram commute
    
    $$
      \xymatrix{
      I\otimes {\cal A}(A,B)
      \ar[rrd]^{=}
      \ar[dd]_{j_{B}\otimes 1}
      &&&&{\cal A}(A,B)\otimes I 
      \ar[dd]^{1\otimes j_{A}}
      \ar[lld]^{=}\\
      &&{\cal A}(A,B)\\
      {\cal A}(B,B)\otimes {\cal A}(A,B)
      \ar[rru]^{M_{ABB}}
      &&&&{\cal A}(A,B)\otimes {\cal A}(A,A)
      \ar[llu]^{M_{AAB}}
      }
   $$
   
   \end{definition}
   
    In general a ${\cal V}$--category is directly analogous to an (ordinary) category enriched over $\mathbf{Set}$ -- 
    if ${\cal V} = \mathbf{Set}$ then these diagrams are the usual category axioms.
    Basically, composition
    of morphisms is replaced by tensoring 
    and the resulting diagrams are required to commute. The next two definitions exhibit this 
    principle and are important since they give 
    us the setting in which to construct a category of ${\cal V}$--categories.
  
   \begin{definition} \label{enriched:funct} For ${\cal V}$--categories 
   ${\cal A}$ and ${\cal B}$, a ${\cal V}$--$functor$ $T:{\cal A}\to{\cal B}$ is a function
    $T:\left| {\cal A} \right| \to \left| {\cal B} \right|$ and a family of 
    morphisms $T_{AB}:{\cal A}(A,B) \to {\cal B}(TA,TB)$ in ${\cal V}$ indexed by 
    pairs $A,B \in \left| {\cal A} \right|$.
    The usual rules for a functor that state $T(f \circ g) = Tf \circ Tg$ 
    and $T1_{A} = 1_{TA}$ become in the enriched setting, respectively, the commuting diagrams
   
   $$
    \xymatrix{
    &{\cal A}(B,C)\otimes {\cal A}(A,B)
    \ar[rr]^{\scriptstyle M}
    \ar[d]^{\scriptstyle T \otimes T}
    &&{\cal A}(A,C)
    \ar[d]^{\scriptstyle T}&\\
    &{\cal B}(TB,TC)\otimes {\cal B}(TA,TB)
    \ar[rr]^{\scriptstyle M}
    &&{\cal B}(TA,TC)
    }
   $$
  and
   $$
    \xymatrix{
    &&{\cal A}(A,A)
    \ar[dd]^{\scriptstyle T_{AA}}\\
    I
    \ar[rru]^{\scriptstyle j_{A}}
    \ar[rrd]_{\scriptstyle j_{TA}}\\
    &&{\cal B}(TA,TA).
    }
   $$
  ${\cal V}$--functors can be composed to form a category called ${\cal V}$--Cat. We will show that this category
  is actually enriched over $\mathbf{Cat}$, the category of (small) categories with cartesian product. 
     \end{definition}
   \clearpage  
  \begin{definition} \label{enr:nat:trans}
  For ${\cal V}$--functors $T,S:{\cal A}\to{\cal B}$ a ${\cal V}$--{\it natural  
   transformation} $\alpha:T \to S:{\cal A} \to {\cal B}$
  is an $\left| {\cal A} \right|$--indexed family of 
  morphisms $\alpha_{A}:I \to {\cal B}(TA,SA)$ satisfying the ${\cal V}$--naturality
  condition expressed by the commutativity of 
  
  $$
    \xymatrix{
    &I \otimes {\cal A}(A,B)
    \ar[rr]^-{\scriptstyle \alpha_{B} \otimes T_{AB}}
    &&{\cal B}(TB,SB) \otimes {\cal B}(TA,TB)
    \ar[rd]^-{\scriptstyle M}
  \\
    {\cal A}(A,B)
    \ar[ru]^{=}
  \ar[rd]_{=}
    &&&&{\cal B}(TA,SB)
  \\
    &{\cal A}(A,B) \otimes I
    \ar[rr]_-{\scriptstyle S_{AB} \otimes \alpha_{A}}
    &&{\cal B}(SA,SB) \otimes {\cal B}(TA,SA)
    \ar[ru]^-{\scriptstyle M}
    }
   $$
  
   \end{definition}
   
   For two ${\cal V}$--functors  $T,S$ to be equal is to say $TA = SA$ for all $A$ 
   and for the ${\cal V}$--natural isomorphism $\alpha$ between them to have 
   components $\alpha_{A} = j_{TA}$. This latter implies equality of the hom--object morphisms: 
   $T_{AB} = S_{AB}$ for all pairs of objects. The implication is seen by combining the second diagram in 
   Definition \ref{V:Cat} with all the diagrams
  in Definitions \ref{enriched:funct} and \ref{enr:nat:trans}.
   
  We want to check that ${\cal V}$--natural transformations can be composed so that  ${\cal V}$--categories, ${\cal V}$--functors
  and ${\cal V}$--natural transformations form a 2--category. First the vertical composite of ${\cal V}$--natural transformations corresponding to the 
  picture
  
  $$
  \xymatrix@R-=16pt{
  &\ar@{=>}[d]^{\alpha}\\
  {\cal A}
  \ar@/^2pc/[rr]^T
  \ar[rr]_>>>>>S
  \ar@/_2pc/[rr]_R
  &\ar@{=>}[d]^{\beta}
  &{\cal B}\\
  &\\
  }
  $$
  has components given by 
  $(\beta \circ \alpha)_{A} = \xymatrix{
                                  I \cong I \otimes I 
                                  \ar[d]_{\beta_{A} \otimes \alpha_{A}}
                                  \\{\cal B}(SA,RA) \otimes {\cal B}(TA,SA)
                                  \ar[d]_M
                                  \\{\cal B}(TA,RA)}$
  
  The reader should check that this composition produces a 
  valid ${\cal V}$--natural transformation and that the composition is associative,
  by using the pentagonal axioms above.
  The identity 2-cells are the identity ${\cal V}$-natural transformations ${\textbf 1}_{Q} : Q \to Q : {\cal B} \to {\cal C}.$
  These are formed from the unit morphisms in ${\cal V}$: ${({\textbf 1}_{Q})}_{B} = j_{QB}$. That this is truly an identity for the 
  vertical composition is easily checked using the second diagram of Definition \ref{V:Cat}.

  In order to define composition of all allowable pasting diagrams in the 2-category, we need to define the 
  composition described by left and right whiskering diagrams and check for independence of 
  choices of order of composition in larger diagrams. The first 
  picture shows a 1-cell (${\cal V}$--functor) following a 2-cell (${\cal V}$--natural transformation).
  These are composed to form a new 2-cell as follows
  
  $$
  
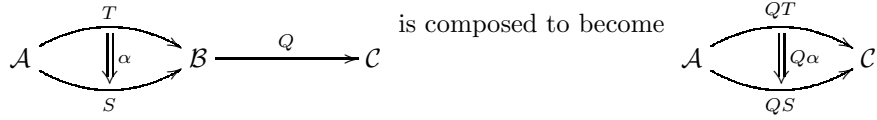
\xymatrix@R-=3pt{
  &\ar@{=>}[dd]^{\alpha}\\
  {\cal A}
  \ar@/^1pc/[rr]^T
  \ar@/_1pc/[rr]_S
  &&{\cal B}
  \ar[rr]^Q
  &&{\cal C}
  \\
  &\\
  } \text{   is composed to become } \xymatrix@R-=3pt{
  &\ar@{=>}[dd]^{Q\alpha}\\
  {\cal A}
  \ar@/^1pc/[rr]^{QT}
  \ar@/_1pc/[rr]_{QS}
  &&{\cal C}
  \\
  &\\
  }
  $$
  
  where $QT$ and $QS$ are given by the usual compositions of their set functions and 
  morphisms in ${\cal V}$, and $Q\alpha$
  has components given by
  $(Q\alpha)_{A} = \xymatrix{
                                  I  
                                  \ar[d]_{\alpha_{A}}
                                  \\{\cal B}(TA,SA)
                                  \ar[d]_{Q_{TA,SA}}
                                  \\{\cal C}(QTA,QSA)}$
                                  
  The second picture shows a 2-cell following a 1-cell. These are composed as follows
  $$
  
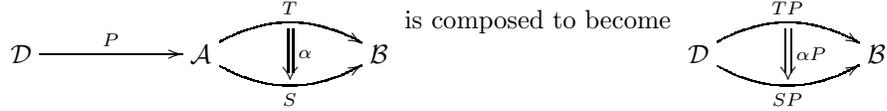
\xymatrix@R-=3pt{
  &&&\ar@{=>}[dd]^{\alpha}\\
  {\cal D}
  \ar[rr]^P
  &&{\cal A}
  \ar@/^1pc/[rr]^T
  \ar@/_1pc/[rr]_S
  &&{\cal B}
  \\
  &&&\\
  } \text{   is composed to become } \xymatrix@R-=3pt{
  &\ar@{=>}[dd]^{\alpha P}\\
  {\cal D}
  \ar@/^1pc/[rr]^{TP}
  \ar@/_1pc/[rr]_{SP}
  &&{\cal B}
  \\
  &\\
  }
  $$
  
  where $\alpha P$ has components given by 
  $(\alpha P)_{D} = \alpha_{PD}$. Again the reader should check that the ${\cal V}$--naturality of $\alpha$ and the
  ${\cal V}$--functoriality of $Q$ imply that the two whisker compositions are ${\cal V}$--natural. 
  What we have developed here are the partial functors of the composition morphism implicit in enriching
  over $\mathbf{Cat}$. The said composition morphism is a functor of two variables. Since $\mathbf{Cat}$ 
  is symmetric, 
  that the partial functors can be combined to make the functor of two variables is implied by the 
  commutativity of a diagram that describes
  the two ways of combining them (see \cite{EK1}). The reader should verify the functoriality of
  the partial functors, while here we will check that they can be combined.
  What needs to be checked is that composing the horizontally adjacent 
  2--cells (shown below) is well--defined and gives whiskering in terms of horizontally composing with an 
  identity 2--cell.
  
  $$
  
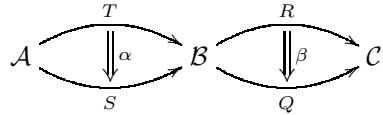
\xymatrix@R-=3pt{
  &\ar@{=>}[dd]^{\alpha}&&\ar@{=>}[dd]^{\beta}\\
  {\cal A}
  \ar@/^1pc/[rr]^T
  \ar@/_1pc/[rr]_S
  &&{\cal B}
  \ar@/^1pc/[rr]^R
  \ar@/_1pc/[rr]_Q
  &&{\cal C}
  \\
  &&&\\
  }
  $$
  
  
  First we need the two ways of composing the above cells using whiskers to be equivalent: 
  $\beta * \alpha = Q\alpha \circ \beta T = \beta S \circ R\alpha$.
  In terms of the above definitions, the following diagram must commute
  
  $$
    \xymatrix{
    &I\otimes I
    \ar[rr]^-{(\beta S)_{A} \otimes (R\alpha)_{A}}
    \ar[d]^{(Q\alpha)_{A} \otimes (\beta T)_{A}}
    &&{\cal C}(RSA,QSA)\otimes {\cal C}(RTA,RSA)
    \ar[d]^{\scriptstyle M}&\\
    &{\cal C}(QTA,QSA)\otimes {\cal C}(RTA,QTA)
    \ar[rr]^{\scriptstyle M}
    &&{\cal C}(RTA,QSA)
    }
   $$
  That this commutes is easily seen since it is just an instance of the diagram in Definition \ref{enr:nat:trans}, specifically 
  with the initial entry being ${\cal B}(TA,SA)$. Associativity of this composition depends on the 
  vertical associativity and on the ${\cal V}$--functoriality of the 1-cells. Since $\mathbf{Cat}$ is 
  symmetric the well--defined nature of the horizontal composition is sufficient to give us all other pasting schemes
  such as, for instance, the
  exchange identity. This states that, in the following picture , 
  $({\beta_2}*{\beta_1})\circ ({\alpha_2}*{\alpha_1}) = ({\beta_2}\circ {\alpha_2})*({\beta_1}\circ {\alpha_1}).$
  
  $$
  
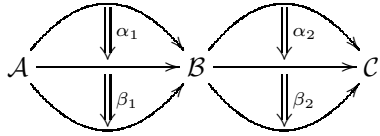
\xymatrix@R-=16pt{
  &\ar@{=>}[d]^{\alpha_1}
  &&\ar@{=>}[d]^{\alpha_2}
  \\
  {\cal A}
  \ar@/^2pc/[rr]
  \ar[rr]
  \ar@/_2pc/[rr]
  &\ar@{=>}[d]^{\beta_1}
  &{\cal B}
  \ar@/^2pc/[rr]
  \ar[rr]
  \ar@/_2pc/[rr]
  &\ar@{=>}[d]^{\beta_2}
  &{\cal C}\\
  &&&&&
  }
  $$
  
  Secondly we need the whiskering to be compatible with horizontal composition with identity 2--cells.
  In other words whiskering 
  a 1--cell Q on the right (or left) of a 2--cell $\alpha: T \to S$  should be the same as horizontally 
  composing ${\textbf 1}_{Q}$ on the respective side of $\alpha$.
  Pictorially for the righthand whiskering:
  $$
  \xymatrix@R-=3pt{
  &\ar@{=>}[dd]^{\alpha}\\
  {\cal A}
  \ar@/^1pc/[rr]^T
  \ar@/_1pc/[rr]_S
  &&{\cal B}
  \ar[rr]^Q
  &&{\cal C}
  \\
  &\\
  } = \xymatrix@R-=3pt{
  &\ar@{=>}[dd]^{\alpha}&&\ar@{=>}[dd]^{{\textbf 1}_{Q}}\\
  {\cal A}
  \ar@/^1pc/[rr]^T
  \ar@/_1pc/[rr]_S
  &&{\cal B}
  \ar@/^1pc/[rr]^Q
  \ar@/_1pc/[rr]_Q
  &&{\cal C}
  \\
  &&&\\
  }
  $$
  To see this equality we need check only one way of composing 
  ${\textbf 1}_{Q} * \alpha$ since we have shown it to be well defined -- i.e. we check that $Q\alpha =  
  {\textbf 1}_{Q} * \alpha = Q\alpha \circ {\textbf 1}_{Q}T$. This is true 
  immediately from the second half of the second diagram
  in Definition \ref{V:Cat}. Now pictorially for the left-hand whiskering:
  
  $$
  \xymatrix@R-=3pt{
  &&&\ar@{=>}[dd]^{\alpha}\\
  {\cal D}
  \ar[rr]^P
  &&{\cal A}
  \ar@/^1pc/[rr]^T
  \ar@/_1pc/[rr]_S
  &&{\cal B}
  \\
  &&&\\
  } = \xymatrix@R-=3pt{
  &\ar@{=>}[dd]^{{\textbf 1}_{P}}&&\ar@{=>}[dd]^{\alpha}\\
  {\cal D}
  \ar@/^1pc/[rr]^P
  \ar@/_1pc/[rr]_P
  &&{\cal A}
  \ar@/^1pc/[rr]^T
  \ar@/_1pc/[rr]_S
  &&{\cal B}
  \\
  &&&\\
  }
  $$
  That $\alpha P = \alpha * {\textbf 1}_{P} = S{\textbf 1}_{P} \circ \alpha P$ is seen from the second diagram of 
  Definition \ref{enriched:funct} for the functor S and the second half of the second diagram in Definition \ref{V:Cat}.
  
  
  Having ascertained that we have a 2--category, it is a good time to review the morphisms between two such things.
  This will make clear what I mean later when I discuss things like a 
  2--functor $\otimes^{(1)}:{\cal V}$--Cat $\times {\cal V}$--Cat $\to {\cal V}$--Cat or a 2--natural transformation
  $c^{(1)}:(-_1\otimes^{(1)}  -_2)\to (-_2\otimes^{(1)} -_1).$ A {\it 2--functor} $F:U\to V$ 
  sends objects to objects, 1--cells to 1--cells, and
  2--cells to 2--cells and preserves all the categorical structure. A {\it 2--natural 
  transformation} $\theta:F\to G:U\to V$ is 
  a function that sends each object $A \in U$ to a 1--cell $\theta_{A}:FA\to GA$ in $V$ in 
  such a way that for each 2--cell
   in $U$ we have that the compositions of the following diagrams are equal in $V$
   
   $$
   \xymatrix@R-=3pt{
   &\ar@{=>}[dd]^{F\alpha}\\
    FA
   \ar@/^1pc/[rr]^{Ff}
   \ar@/_1pc/[rr]_{Fg}
   &&FB
   \ar[rr]^{\theta_B}
   &&GB
   \\
  &\\
  }
  =
  \xymatrix@R-=3pt{
  &&&\ar@{=>}[dd]^{G\alpha}\\
  FA
  \ar[rr]^{\theta_A}
  &&GA
  \ar@/^1pc/[rr]^{Gf}
  \ar@/_1pc/[rr]_{Gg}
  &&GB
  \\
  &&&\\
  }
  $$
  
  Furthermore for two 2--natural transformations a {\it modification} $\mu:\theta\to \phi:F\to G:U\to V$ is a 
  function that sends each object $A \in U$ to a 2--cell $\mu_A:\theta_{A}\to \phi_A:FA\to GA$ in such a way that 
  for each 2--cell
   in $U$ we have that the compositions of the following diagrams are equal in $V$
   
   $$
   \xymatrix@R-=3pt{
   &\ar@{=>}[dd]^{F\alpha}&&\ar@{=>}[dd]^{\mu_B}\\
   FA
   \ar@/^1pc/[rr]^{Ff}
   \ar@/_1pc/[rr]_{Fg}
   &&FB
   \ar@/^1pc/[rr]^{\theta_B}
   \ar@/_1pc/[rr]_{\phi_B}
   &&GB
   \\
   &&&\\
   }
  =
  \xymatrix@R-=3pt{
  &\ar@{=>}[dd]^{\mu_A}&&\ar@{=>}[dd]^{G\alpha}\\
  FA
  \ar@/^1pc/[rr]^{\theta_A}
  \ar@/_1pc/[rr]_{\phi_A}
  &&GA
  \ar@/^1pc/[rr]^{Gf}
  \ar@/_1pc/[rr]_{Gg}
  &&GB
  \\
  &&&\\
  }
  $$
  Taking modifications as 3--cells, 2--natural transformations as 2--cells, 2--functors as 1--cells, and 2--categories 
  as objects gives us a 3--category called 2--Cat.

  
  \clearpage
  \newpage
    \MySection{Categories Enriched over a Braided Monoidal Category}
  
  \begin{definition} A {\it braiding} for a monoidal category ${\cal V}$ is a family of natural  
  isomorphisms $c_{XY}: X \otimes Y \to Y \otimes X$
  such that the following diagrams commute. They are drawn next to their underlying braids.
  \begin{enumerate}
  \item
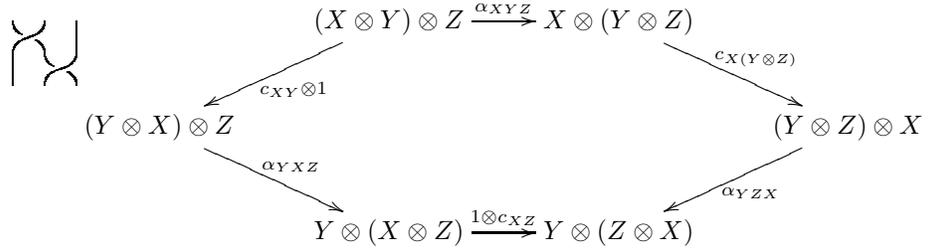
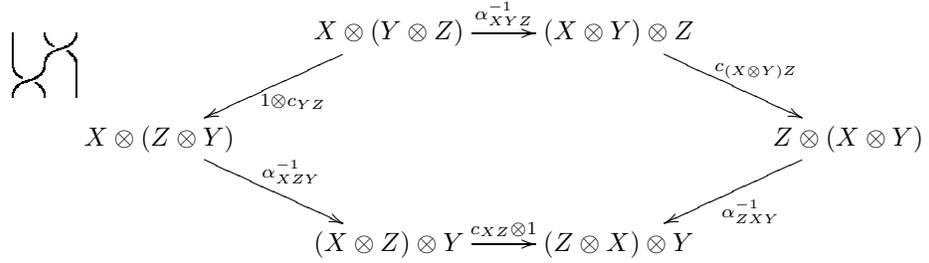
  $$
  \xy 0;/r1pc/:
  ,{\vtwist}+(2,1),{\xcapv-@(0)}
  +(-2,0),{\xcapv-@(0)}+(1,1),{\vtwist}
  \endxy
  \xymatrix{
  &(X \otimes Y) \otimes Z
  \ar[dl]^{c_{XY} \otimes 1}
  \ar[r]^{\alpha_{XYZ}}
  &X \otimes (Y \otimes Z)
  \ar[dr]^{c_{X(Y \otimes Z)}}
  \\(Y \otimes X) \otimes Z
  \ar[dr]^{\alpha_{YXZ}}
  &&&(Y \otimes Z) \otimes X
  \ar[dl]^{\alpha_{YZX}}
  \\&Y \otimes (X \otimes Z)
  \ar[r]^{1\otimes c_{XZ}}
  &Y \otimes (Z \otimes X)
  }
  $$

  \item
  $$
  \xy 0;/r1pc/:
  ,{\xcapv-@(0)}+(1,1),{\vtwist}+(-1,0)
  ,{\vtwist}+(2,1),{\xcapv-@(0)}
  \endxy
  \xymatrix{
  &X \otimes (Y \otimes Z)
  \ar[dl]^{1 \otimes c_{YZ}}
  \ar[r]^{\alpha_{XYZ}^{-1}}
  &(X \otimes Y) \otimes Z
  \ar[dr]^{c_{(X \otimes Y)Z}}
  \\X \otimes (Z \otimes Y)
  \ar[dr]^{\alpha_{XZY}^{-1}}
  &&&Z \otimes (X \otimes Y)
  \ar[dl]^{\alpha_{ZXY}^{-1}}
  \\&(X \otimes Z) \otimes Y
  \ar[r]^{c_{XZ} \otimes 1}
  &(Z \otimes X) \otimes Y
  }
  $$
  \end{enumerate}
  
  A braided category is a monoidal category with a chosen braiding.
  \end{definition}
  Joyal and Street 
  proved the coherence theorem
  for braided categories in \cite{JS}, the immediate corollary of which is that in a free 
  braided category generated by a set of
  objects, a diagram commutes if and only if all legs having the same source and target have the same underlying braid.

  \begin{definition}
  A {\it symmetry} is a braiding such that the following diagram commutes
  $$
  
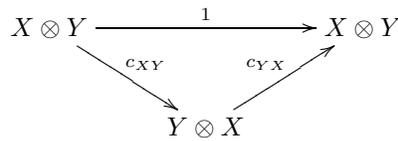
\xymatrix{
  X \otimes Y
  \ar[rr]^1
  \ar[dr]^{c_{XY}}
  &&X \otimes Y
  \\
  &Y \otimes X
  \ar[ur]^{c_{YX}}
  }
  $$
  In other words $c_{XY}^{-1} = c_{YX}$. A symmetric category is a monoidal category with a chosen symmetry.
  \end{definition}
  \clearpage
  If ${\cal V}$ is braided then we can define additional structure on ${\cal V}$--Cat. First
  there is an opposite of a ${\cal V}$--category which has 
  $\left|{\cal A}^{op} \right| = \left|{\cal A}\right|$ and ${\cal A}^{op}(A,A') = {\cal A}(A',A).$ The
  composition morphisms are given by
  $$
    \xymatrix{
    {\cal A}^{op}(A',A'')\otimes {\cal A}^{op}(A,A')
    \ar@{=}[d]\\
    {\cal A}(A'',A')\otimes {\cal A}(A',A)
    \ar[d]_{c_{{\cal A}(A'',A')\otimes {\cal A}(A',A)}}\\
    {\cal A}(A',A)\otimes {\cal A}(A'',A')
    \ar[d]_{M_{AA'A''}}\\
    {\cal A}(A'',A)
    \ar@{=}[d]\\
    {\cal A}^{op}(A,A'')
    }
  $$
  It is clear from this that $({\cal A}^{op})^{op}\ne {\cal A}.$ 
  The pentagon diagram for the composition morphisms commutes since the braids underlying 
  its legs are the two sides of the braid relation, also known as the Yang-Baxter equation.
  The second structure is a
  product for ${\cal V}$--Cat, that is, a 2-functor 
  $$\otimes^{(1)} : {\cal V} \text{--Cat} \times {\cal V} \text{--Cat} \to {\cal V} \text{--Cat}.$$
  I will always denote the product(s) in ${\cal V}$--Cat
  with a superscript in parentheses that corresponds to the level of enrichment of the components of their domain.
  The product(s) in ${\cal V}$ 
  should logically then have a superscript (0) but I have suppressed this for brevity and to agree with my sources. 
  The product of two ${\cal V}$--categories ${\cal A}$ and ${\cal B}$ has 
  $\left|{\cal A} \otimes^{(1)} {\cal B}\right| = \left|{\cal A}\right| \times \left|{\cal B}\right|$ and 
  $({\cal A} \otimes^{(1)} {\cal B})((A,B),(A',B')) = {\cal A}(A,A') \otimes {\cal B}(B,B').$ 
  
  The unit morphisms for the product ${\cal V}$--categories are the composites
  $$\xymatrix{I \cong I \otimes I \ar[r]_<<<<<{j_{A} \otimes j_{B}} & {\cal A}(A,A)\otimes {\cal B}(B,B) }$$
  
  
  The composition morphisms 
  \begin{small}
  $$M_{(A,B)(A',B')(A'',B'')} : ({\cal A}\otimes^{(1)} {\cal B})((A',B'),(A'',B''))\otimes ({\cal A}\otimes^{(1)} {\cal B})((A,B),(A',B'))\to ({\cal A}\otimes^{(1)} {\cal B})((A,
  B),(A'',B''))$$
  \end{small}
  may be given by
  $$
  \xymatrix{
  ({\cal A}\otimes^{(1)} {\cal B})((A',B'),(A'',B''))\otimes ({\cal A}\otimes^{(1)} {\cal B})((A,B),(A',B'))
  \ar@{=}[d]\\
  ({\cal A}(A',A'')\otimes {\cal B}(B',B''))\otimes ({\cal A}(A,A')\otimes {\cal B}(B,B'))
  \ar[d]_{(1 \otimes \alpha^{-1}) \circ \alpha}\\
  {\cal A}(A',A'')\otimes (({\cal B}(B',B'')\otimes {\cal A}(A,A'))\otimes {\cal B}(B,B'))
  \ar[d]_{1 \otimes (c_{{\cal B}(B',B''){\cal A}(A,A')} \otimes 1)}\\
  {\cal A}(A',A'')\otimes (({\cal A}(A,A'))\otimes {\cal B}(B',B'')) \otimes {\cal B}(B,B'))
  \ar[d]_{\alpha^{-1} \circ (1 \otimes \alpha)}\\
  ({\cal A}(A',A'')\otimes {\cal A}(A,A'))\otimes ({\cal B}(B',B'')) \otimes {\cal B}(B,B'))
  \ar[d]_{M_{AA'A''}\otimes M_{BB'B''}}\\
  ({\cal A}(A,A'')\otimes {\cal B}(B,B''))
  \ar@{=}[d]\\
  ({\cal A}\otimes^{(1)} {\cal B})((A,B),(A'',B''))
  }
  $$
  
  That $({\cal A}\otimes^{(1)} {\cal B})^{op} \ne {\cal A}^{op}\otimes^{(1)} {\cal B}^{op}$ follows from the 
  following braid inequality:
  $$
  \begin{turn}{180}
    \xy 0;/r1pc/:
    ,{\vtwist}+(2,1),{\vtwist}+(-2,0)
    ,{\xcapv-@(0)}+(1,1),{\vtwist}+(2,1),{\xcapv-@(0)}+(-3,0)
    \endxy
    \xymatrix{&\ne}
    \xy 0;/r1pc/:
    ,{\xcapv-@(0)}+(1,1),{\vtwist}+(2,1),{\xcapv-@(0)}+(-3,0)
    ,{\xcapv-@(0)}+(1,1),{\vtwist}+(2,1),{\xcapv-@(0)}+(-3,0)
    ,{\vtwist}+(2,1),{\vtwist}+(-2,0)
    ,{\xcapv-@(0)}+(1,1),{\vtwist}+(2,1),{\xcapv-@(0)}+(-3,0)
    \endxy
   \end{turn} 
  $$
  
  Now consider more carefully the morphisms of ${\cal V}$ that make up the compostion morphism for a product
  enriched category, especially those
  that accomplish the interchange of the interior hom--objects.
  In the symmetric case, any other combination of instances of $\alpha$ and $c$ 
  with the same domain and range would be equal, due to coherence. 
  In the merely braided case, there at first seems to be a much larger range of available choices. There is a canonical epimorphism 
  $\sigma : B_n \to S_n$ of the braid group on n strands onto the permutation group. The permutation given by $\sigma$ is that 
  given by the strands of the braid on the $n$ original positions. For instance on a canonical generator of $B_n$, 
  $\sigma_i$, we have $\sigma(\sigma_i) = (i\text{ }\text{ }i+1)$.
  Candidates for multiplication would seem to be those defined using any braid $b \in B_4$ 
  such that $\sigma(b) = (2 3)$. 
  It is clear that the composition morphism would be defined as above, with a series of instances of $\alpha$ and $c$
  such that the underlying braid is $b$, followed in turn  by $M_{AA'A''}\otimes M_{BB'B''}$ in order to complete the 
  composition. That $M_{AA'A''}\otimes M_{BB'B''}$ will have the correct domain on which to operate is guaranteed by the 
  permutation condition 
  on $b$.
  \clearpage
  The unit axioms hold due to the naturality of compositions of $\alpha$ and $c$ and the unit axioms obeyed by 
  ${\cal A}$ and ${\cal B}.$
  The remaining things to be checked are associativity of composition and functoriality of the associator. 
  
  For the associativity axiom to hold the following diagram must 
    commute, where the initial bullet represents
     $$
    [({\cal A}\otimes^{(1)} {\cal B})((A'',B''),(A''',B'''))\otimes ({\cal A}\otimes^{(1)} {\cal B})((A',B'),(A'',B''))]\otimes ({\cal A}\otimes^{(1)} {\cal B})((A,B),(A',B'))
    $$
    and the last bullet represents
    $[{\cal A}\otimes^{(1)} {\cal B}]((A,B),(A''',B''')).$
    
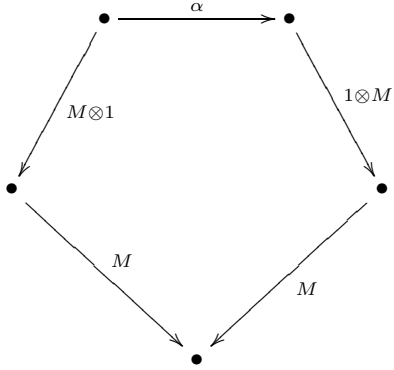
    $$
      \xymatrix{
      &\bullet
      \ar[rr]^{ \alpha}
      \ar[ddl]^{ M \otimes 1}
      &&\bullet
      \ar[ddr]^{ 1 \otimes M}&\\\\
      \bullet
      \ar[ddrr]^{ M}
      &&&&\bullet
      \ar[ddll]^{ M}
      \\\\&&\bullet
      }$$

    In ${\cal V}$ let 
    $X={\cal A}(A,A')$, $X'={\cal A}(A',A'')$, $X''={\cal A}(A'',A''')$, $Y={\cal B}(B,B')$, $Y'={\cal B}(B',B'')$ and $Y''={\cal B}(B'',B''').$
    The exterior of the following expanded diagram (where I leave out some parentheses 
    for clarity and denote various composites
  of $\alpha$ and $c$ by unlabeled arrows) is required to commute. 
    
    \noindent
              	          \begin{center}
    	          \resizebox{6.5in}{!}
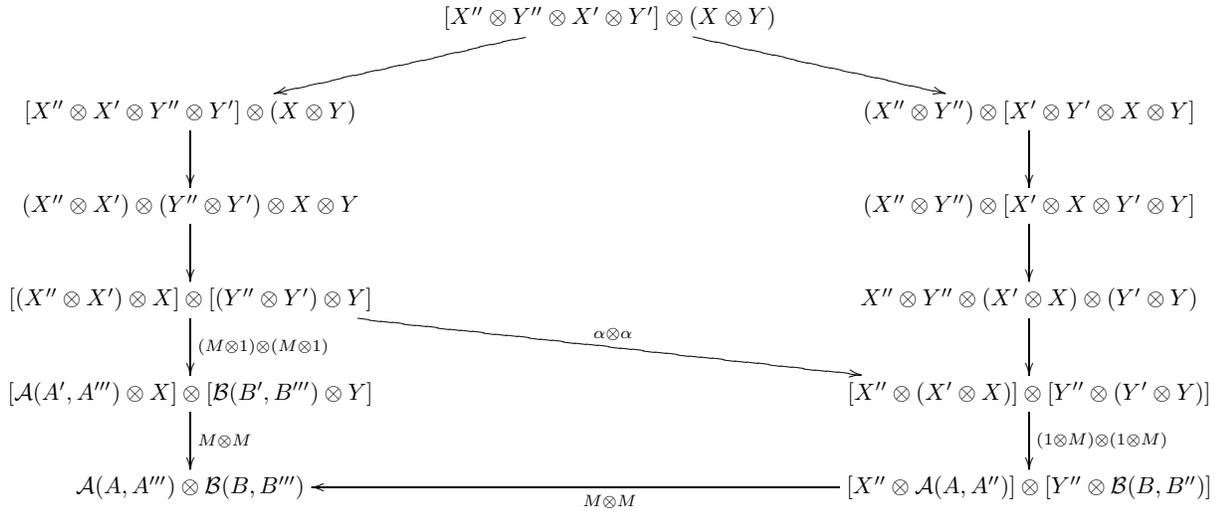
{
    $$
    \xymatrix{
    &[X''\otimes Y''\otimes X'\otimes Y']\otimes (X\otimes Y)
    \ar[dr]
    \ar[dl]\\
    [X''\otimes X'\otimes Y''\otimes Y']\otimes (X\otimes Y)
    \ar[d]
    &&(X''\otimes Y'')\otimes [X'\otimes Y'\otimes X\otimes Y]
    \ar[d]\\
    (X''\otimes X')\otimes (Y''\otimes Y')\otimes X\otimes Y
    \ar[d]
    &&(X''\otimes Y'')\otimes [X'\otimes X\otimes Y'\otimes Y]
    \ar[d]
    \\
    [(X''\otimes X')\otimes X]\otimes [(Y''\otimes Y')\otimes Y]
        \ar[drr]^{\alpha \otimes \alpha}
    \ar[d]^{(M\otimes 1)\otimes (M \otimes 1)}
    &&X''\otimes Y''\otimes (X'\otimes X)\otimes (Y'\otimes Y)
    \ar[d]\\
    [{\cal A}(A',A''')\otimes X]\otimes [{\cal B}(B',B''')\otimes Y]
    \ar[d]^{M\otimes M}
    &&[X''\otimes (X'\otimes X)]\otimes [Y''\otimes (Y'\otimes Y)]
    \ar[d]^{(1\otimes M)\otimes (1 \otimes M)}\\
    {\cal A}(A,A''')\otimes {\cal B}(B,B''')
    &&[X''\otimes {\cal A}(A,A'')]\otimes [Y''\otimes {\cal B}(B,B'')]
    \ar[ll]^{M\otimes M}
    }
    $$
    }
    \end{center}
    
    The bottom region commutes by the associativity axioms for ${\cal A}$ and ${\cal B}.$
    \clearpage
    We are left needing to show that the underlying braids are equal
    for the two legs of the upper region.
    Again these basic nodes must be present regardless of the choice of braid by which the composition morphism is defined.  
    Notice that  
    the right and left legs have the following underlying braids in $B_6$
    for some examples of various choices of $b$ in $B_4.$ We call the two derived braids in $B_6$ $Lb$ and $Rb$
    respectively. $Lb$ is algorithmically described as a copy of $b$ on the first 4 strands followed by 
    a copy of $b$ on the 4 ``strands'' that result from pairing as the edges of two ribbons strands 1 and 2, and strands 
    3 and 4, along with the remaining two strands 5 and 6. $Rb$ is similarly described, but the initial
    copy of $b$ is on the last 4 strands, and the ribbon edge pairing is on the pairs 4,5 and 5,6.
    The first example for $b$ is the one used in the original definition given above.

     $$
     b_{(1)}=
     \begin{turn}{180}
     \xy 0;/r1pc/:
     ,{\xcapv-@(0)}+(1,1),{\vtwist}+(2,1),{\xcapv-@(0)}
     \endxy
     \end{turn}
     \text{ Associativity follows from: }
     \begin{turn}{180}
     \xy 0;/r1pc/:
        ,{\xcapv-@(0)}+(1,1),{\xcapv-@(0)}+(1,1),{\vtwist}+(2,1),{\xcapv-@(0)}+(1,1),{\xcapv-@(0)}+(-5,0)
        ,{\xcapv-@(0)}+(1,1),{\vtwist}+(2,1),{\xcapv-@(0)}+(1,1),{\xcapv-@(0)}+(1,1),{\xcapv-@(0)}+(-5,0)
        ,{\xcapv-@(0)}+(1,1),{\xcapv-@(0)}+(1,1),{\xcapv-@(0)}+(1,1),{\vtwist}+(2,1),{\xcapv-@(0)}
     \endxy
     \end{turn}
     \xymatrix{=}
     \begin{turn}{180}
     \xy 0;/r1pc/:
        ,{\xcapv-@(0)}+(1,1),{\xcapv-@(0)}+(1,1),{\vtwist}+(2,1),{\xcapv-@(0)}+(1,1),{\xcapv-@(0)}+(-5,0)
        ,{\xcapv-@(0)}+(1,1),{\xcapv-@(0)}+(1,1),{\xcapv-@(0)}+(1,1),{\vtwist}+(2,1),{\xcapv-@(0)}+(-5,0)
        ,{\xcapv-@(0)}+(1,1),{\vtwist}+(2,1),{\xcapv-@(0)}+(1,1),{\xcapv-@(0)}+(1,1),{\xcapv-@(0)}
     \endxy
     \end{turn}
     $$
     
     $$
     b_{(2)}=
     \begin{turn}{180}
     \xy 0;/r1pc/:
     ,{\xcapv-@(0)}+(1,1),{\vtwist}+(2,1),{\xcapv-@(0)}+(-3,0)
     ,{\vtwist}+(2,1),{\xcapv-@(0)}+(1,1),{\xcapv-@(0)}+(-3,0)
     ,{\vtwist}+(2,1),{\xcapv-@(0)}+(1,1),{\xcapv-@(0)}
     \endxy
     \end{turn}
     \text{ Associativity does not follow since: }
     \begin{turn}{180}
     \xy 0;/r1pc/:
        ,{\xcapv-@(0)}+(1,1),{\xcapv-@(0)}+(1,1),{\vtwist}+(2,1),{\xcapv-@(0)}+(1,1),{\xcapv-@(0)}+(-5,0)
        ,{\xcapv-@(0)}+(1,1),{\vtwist}+(2,1),{\xcapv-@(0)}+(1,1),{\xcapv-@(0)}+(1,1),{\xcapv-@(0)}+(-5,0)
        ,{\vtwist}+(2,1),{\xcapv-@(0)}+(1,1),{\xcapv-@(0)}+(1,1),{\xcapv-@(0)}+(1,1),{\xcapv-@(0)}+(-5,0)
        ,{\vtwist}+(2,1),{\xcapv-@(0)}+(1,1),{\xcapv-@(0)}+(1,1),{\xcapv-@(0)}+(1,1),{\xcapv-@(0)}+(-5,0)
        ,{\xcapv-@(0)}+(1,1),{\xcapv-@(0)}+(1,1),{\xcapv-@(0)}+(1,1),{\vtwist}+(2,1),{\xcapv-@(0)}+(-5,0)
        ,{\xcapv-@(0)}+(1,1),{\xcapv-@(0)}+(1,1),{\vtwist}+(2,1),{\xcapv-@(0)}+(1,1),{\xcapv-@(0)}+(-5,0)
        ,{\xcapv-@(0)}+(1,1),{\xcapv-@(0)}+(1,1),{\vtwist}+(2,1),{\xcapv-@(0)}+(1,1),{\xcapv-@(0)}+(-5,0)
     \endxy
     \end{turn}
     \xymatrix{\ne}
     \begin{turn}{180}
     \xy 0;/r1pc/:
        ,{\xcapv-@(0)}+(1,1),{\xcapv-@(0)}+(1,1),{\vtwist}+(2,1),{\xcapv-@(0)}+(1,1),{\xcapv-@(0)}+(-5,0)
        ,{\xcapv-@(0)}+(1,1),{\vtwist}+(2,1),{\vtwist}+(2,1),{\xcapv-@(0)}+(-5,0)
        ,{\vtwist}+(2,1),{\vtwist}+(2,1),{\xcapv-@(0)}+(1,1),{\xcapv-@(0)}+(-5,0)
        ,{\xcapv-@(0)}+(1,1),{\vtwist}+(2,1),{\xcapv-@(0)}+(1,1),{\xcapv-@(0)}+(1,1),{\xcapv-@(0)}+(-5,0)
        ,{\xcapv-@(0)}+(1,1),{\vtwist}+(2,1),{\xcapv-@(0)}+(1,1),{\xcapv-@(0)}+(1,1),{\xcapv-@(0)}+(-5,0)
        ,{\vtwist}+(2,1),{\vtwist}+(2,1),{\xcapv-@(0)}+(1,1),{\xcapv-@(0)}+(-5,0)
        ,{\xcapv-@(0)}+(1,1),{\vtwist}+(2,1),{\xcapv-@(0)}+(1,1),{\xcapv-@(0)}+(1,1),{\xcapv-@(0)}+(-5,0)
        ,{\xcapv-@(0)}+(1,1),{\vtwist}+(2,1),{\xcapv-@(0)}+(1,1),{\xcapv-@(0)}+(1,1),{\xcapv-@(0)}+(-5,0)
        ,{\vtwist}+(2,1),{\xcapv-@(0)}+(1,1),{\xcapv-@(0)}+(1,1),{\xcapv-@(0)}+(1,1),{\xcapv-@(0)}+(-5,0)
        ,{\vtwist}+(2,1),{\xcapv-@(0)}+(1,1),{\xcapv-@(0)}+(1,1),{\xcapv-@(0)}+(1,1),{\xcapv-@(0)}+(-5,0)
     \endxy
     \end{turn}
     $$
     
     $$
     b_{(3)}=
     \begin{turn}{180}
     \xy 0;/r1pc/:
     ,{\xcapv-@(0)}+(1,1),{\xcapv-@(0)}+(1,1),{\vtwist}+(-2,0)
     ,{\xcapv-@(0)}+(1,1),{\xcapv-@(0)}+(1,1),{\vtwist}+(-2,0)
     ,{\xcapv-@(0)}+(1,1),{\vtwist}+(2,1),{\xcapv-@(0)}
     \endxy
     \end{turn}
     \text{ Associativity follows from: }
     \begin{turn}{180}
     \xy 0;/r1pc/:
        ,{\xcapv-@(0)}+(1,1),{\xcapv-@(0)}+(1,1),{\xcapv-@(0)}+(1,1),{\vtwist}+(2,1),{\xcapv-@(0)}+(-5,0)
        ,{\xcapv-@(0)}+(1,1),{\xcapv-@(0)}+(1,1),{\xcapv-@(0)}+(1,1),{\xcapv-@(0)}+(1,1),{\vtwist}+(-4,0)
        ,{\xcapv-@(0)}+(1,1),{\xcapv-@(0)}+(1,1),{\xcapv-@(0)}+(1,1),{\xcapv-@(0)}+(1,1),{\vtwist}+(-4,0)
        ,{\xcapv-@(0)}+(1,1),{\xcapv-@(0)}+(1,1),{\xcapv-@(0)}+(1,1),{\vtwist}+(2,1),{\xcapv-@(0)}+(-5,0)
        ,{\xcapv-@(0)}+(1,1),{\xcapv-@(0)}+(1,1),{\vtwist}+(2,1),{\xcapv-@(0)}+(1,1),{\xcapv-@(0)}+(-5,0)
        ,{\xcapv-@(0)}+(1,1),{\vtwist}+(2,1),{\xcapv-@(0)}+(1,1),{\xcapv-@(0)}+(1,1),{\xcapv-@(0)}+(-5,0)
        ,{\xcapv-@(0)}+(1,1),{\xcapv-@(0)}+(1,1),{\xcapv-@(0)}+(1,1),{\xcapv-@(0)}+(1,1),{\vtwist}+(-4,0)
        ,{\xcapv-@(0)}+(1,1),{\xcapv-@(0)}+(1,1),{\xcapv-@(0)}+(1,1),{\xcapv-@(0)}+(1,1),{\vtwist}+(-4,0)
        ,{\xcapv-@(0)}+(1,1),{\xcapv-@(0)}+(1,1),{\xcapv-@(0)}+(1,1),{\vtwist}+(2,1),{\xcapv-@(0)}+(-5,0)
     \endxy
     \end{turn}
     \xymatrix{=}
     \begin{turn}{180}
     \xy 0;/r1pc/:
        ,{\xcapv-@(0)}+(1,1),{\xcapv-@(0)}+(1,1),{\xcapv-@(0)}+(1,1),{\xcapv-@(0)}+(1,1),{\vtwist}+(-4,0)
        ,{\xcapv-@(0)}+(1,1),{\xcapv-@(0)}+(1,1),{\xcapv-@(0)}+(1,1),{\vtwist}+(2,1),{\xcapv-@(0)}+(-5,0)
        ,{\xcapv-@(0)}+(1,1),{\xcapv-@(0)}+(1,1),{\xcapv-@(0)}+(1,1),{\vtwist}+(2,1),{\xcapv-@(0)}+(-5,0)
        ,{\xcapv-@(0)}+(1,1),{\xcapv-@(0)}+(1,1),{\vtwist}+(2,1),{\vtwist}+(-4,0)
        ,{\xcapv-@(0)}+(1,1),{\xcapv-@(0)}+(1,1),{\xcapv-@(0)}+(1,1),{\vtwist}+(2,1),{\xcapv-@(0)}+(-5,0)
        ,{\xcapv-@(0)}+(1,1),{\xcapv-@(0)}+(1,1),{\vtwist}+(2,1),{\xcapv-@(0)}+(1,1),{\xcapv-@(0)}+(-5,0)
        ,{\xcapv-@(0)}+(1,1),{\xcapv-@(0)}+(1,1),{\vtwist}+(2,1),{\xcapv-@(0)}+(1,1),{\xcapv-@(0)}+(-5,0)
        ,{\xcapv-@(0)}+(1,1),{\vtwist}+(2,1),{\xcapv-@(0)}+(1,1),{\xcapv-@(0)}+(1,1),{\xcapv-@(0)}
     \endxy
     \end{turn}
     $$
     \clearpage
     $$
     b_{(4)}=
     \begin{turn}{180}
     \xy 0;/r1pc/:
     ,{\xcapv-@(0)}+(1,1),{\vtwist}+(2,1),{\xcapv-@(0)}+(-3,0)
     ,{\xcapv-@(0)}+(1,1),{\vtwist}+(2,1),{\xcapv-@(0)}+(-3,0)
     ,{\xcapv-@(0)}+(1,1),{\vtwist}+(2,1),{\xcapv-@(0)}+(-3,0)
     \endxy
     \end{turn}
     \text{ Associativity does not follow since: }
     \begin{turn}{180}
     \xy 0;/r1pc/:
        ,{\xcapv-@(0)}+(1,1),{\xcapv-@(0)}+(1,1),{\vtwist}+(2,1),{\xcapv-@(0)}+(1,1),{\xcapv-@(0)}+(-5,0)
        ,{\xcapv-@(0)}+(1,1),{\vtwist}+(2,1),{\xcapv-@(0)}+(1,1),{\xcapv-@(0)}+(1,1),{\xcapv-@(0)}+(-5,0)
        ,{\xcapv-@(0)}+(1,1),{\vtwist}+(2,1),{\xcapv-@(0)}+(1,1),{\xcapv-@(0)}+(1,1),{\xcapv-@(0)}+(-5,0)
        ,{\xcapv-@(0)}+(1,1),{\xcapv-@(0)}+(1,1),{\vtwist}+(2,1),{\xcapv-@(0)}+(1,1),{\xcapv-@(0)}+(-5,0)
        ,{\xcapv-@(0)}+(1,1),{\xcapv-@(0)}+(1,1),{\vtwist}+(2,1),{\xcapv-@(0)}+(1,1),{\xcapv-@(0)}+(-5,0)
        ,{\xcapv-@(0)}+(1,1),{\vtwist}+(2,1),{\xcapv-@(0)}+(1,1),{\xcapv-@(0)}+(1,1),{\xcapv-@(0)}+(-5,0)
        ,{\xcapv-@(0)}+(1,1),{\xcapv-@(0)}+(1,1),{\xcapv-@(0)}+(1,1),{\vtwist}+(2,1),{\xcapv-@(0)}+(-5,0)
        ,{\xcapv-@(0)}+(1,1),{\xcapv-@(0)}+(1,1),{\xcapv-@(0)}+(1,1),{\vtwist}+(2,1),{\xcapv-@(0)}+(-5,0)
        ,{\xcapv-@(0)}+(1,1),{\xcapv-@(0)}+(1,1),{\xcapv-@(0)}+(1,1),{\vtwist}+(2,1),{\xcapv-@(0)}+(-5,0)
     \endxy
     \end{turn}
     \xymatrix{\ne}
     \begin{turn}{180}
     \xy 0;/r1pc/:
        ,{\xcapv-@(0)}+(1,1),{\xcapv-@(0)}+(1,1),{\vtwist}+(2,1),{\xcapv-@(0)}+(1,1),{\xcapv-@(0)}+(-5,0)
        ,{\xcapv-@(0)}+(1,1),{\xcapv-@(0)}+(1,1),{\xcapv-@(0)}+(1,1),{\vtwist}+(2,1),{\xcapv-@(0)}+(-5,0)
        ,{\xcapv-@(0)}+(1,1),{\xcapv-@(0)}+(1,1),{\xcapv-@(0)}+(1,1),{\vtwist}+(2,1),{\xcapv-@(0)}+(-5,0)
        ,{\xcapv-@(0)}+(1,1),{\xcapv-@(0)}+(1,1),{\vtwist}+(2,1),{\xcapv-@(0)}+(1,1),{\xcapv-@(0)}+(-5,0)
        ,{\xcapv-@(0)}+(1,1),{\xcapv-@(0)}+(1,1),{\vtwist}+(2,1),{\xcapv-@(0)}+(1,1),{\xcapv-@(0)}+(-5,0)
        ,{\xcapv-@(0)}+(1,1),{\xcapv-@(0)}+(1,1),{\xcapv-@(0)}+(1,1),{\vtwist}+(2,1),{\xcapv-@(0)}+(-5,0)
        ,{\xcapv-@(0)}+(1,1),{\vtwist}+(2,1),{\xcapv-@(0)}+(1,1),{\xcapv-@(0)}+(1,1),{\xcapv-@(0)}+(-5,0)
        ,{\xcapv-@(0)}+(1,1),{\vtwist}+(2,1),{\xcapv-@(0)}+(1,1),{\xcapv-@(0)}+(1,1),{\xcapv-@(0)}+(-5,0)
        ,{\xcapv-@(0)}+(1,1),{\vtwist}+(2,1),{\xcapv-@(0)}+(1,1),{\xcapv-@(0)}+(1,1),{\xcapv-@(0)}+(-5,0)
     \endxy
     \end{turn}
     $$

Before turning to check on functoriality of the associator, we note that $b_{(3)}$ is the braid underlying the
composition morphism of the product category $({\cal A}^{op})^{op}\otimes^{(1)} {\cal B}$ where
the product is defined using $b_{(1)}.$ This provides the
hint that the two derived braids in $B_6$ that we are comparing above are equal because of the fact that
the opposite of a ${\cal V}$--category is a valid ${\cal V}$--category. In fact we can describe
sufficient conditions for $Lb$ to be equivalent to $Rb$ by describing the braids $b$ that underlie
the composition morphism of a product category given generally by 
$((({\cal A}^{op})^{...op}\otimes^{(1)} ({\cal B}^{op})^{...op})^{op})^{...op}$ where the number of
$op$  exponents is arbitrary in each position. Those braids are alternately described as lying in
$H\sigma_2K \subset B_4$ where $H$ is the cyclic subgroup generated by the braid $\sigma_2\sigma_1\sigma_3\sigma_2$
and $K$ is the subgroup generated by the two generators $\{\sigma_1, \sigma_3\}.$ 
The latter subgroup $K$ is isomorphic to $Z\times Z.$ The first coordinate corresponds to the number of 
$op$ exponents on ${\cal A}$ and the second component to the number of $op$ exponents on ${\cal B}.$ The
power of the element of $H$ corresponds to the 
number of $op$ exponents on the product of the two enriched categories, that is, the number of 
$op$ exponents outside the parentheses. That $b\in H\sigma_2K$ implies $Lb=Rb$ follows from the 
fact that the composition morphisms belonging to the opposite of a ${\cal V}$--category obey the 
pentagon axiom. An exercise of some value is to check consistency of the definitions by constructing
an inductive proof of the implication based on braid group generators.
This is not a necessary condition,
but it may be when the additional requirement that $\sigma(b) = (2 3)$ is added. More work needs to be done to 
determine the necessary conditions and to study the structure and properties of the braids that meet these
conditions.

  Functoriality of the associator is necessary because
  here we need a 2--natural transformation $\alpha^{(1)}$. This means we have a family of ${\cal V}$--functors
  indexed by triples of ${\cal V}$--categories. On objects $\alpha^{(1)}_{{\cal A}{\cal B}{\cal C}}((A,B),C) = (A,(B,C)).$
  In order to guarantee that $\alpha^{(1)}$ obey the coherence pentagon for hom--object morphisms, 
  we define it to be {\it based upon}
  $\alpha$ in ${\cal V}.$ This means precisely that:
  \begin{small}
  $$\alpha^{(1)}_{{\cal A}{\cal B}{\cal C}_{((A,B),C)((A',B'),C')}}: [({\cal A}\otimes^{(1)} {\cal B})\otimes^{(1)} {\cal C}](((A,B),C)((A',B'),C')) \to [{\cal A}\otimes^{(1)} ({\cal B}\otimes^{(1)} {\cal C})]((A,(B,C))(A',(B',C')))$$
  \end{small}
  is equal to
  \begin{small}
  $$\alpha_{{\cal A}(A,A'){\cal B}(B,B'){\cal C}(C,C')}:({\cal A}(A,A')\otimes {\cal B}(B,B'))\otimes {\cal C}(C,C')\to {\cal A}(A,A') \otimes ({\cal B}(B,B')\otimes {\cal C}(C,C')).$$
  \end{small}
  \clearpage
  This definition guarantees that the $\alpha^{(1)}$ pentagons for objects and for hom--objects commute: 
  the first trivially and the second by the fact that the
  $\alpha$ pentagon commutes in ${\cal V}.$
  We must also check for ${\cal V}$--functoriality. The unit axioms are trivial -- we consider the more interesting 
  associativity of hom--object morphisms property. The following diagram must commute, where the first bullet represents
  $$[({\cal A}\otimes^{(1)} {\cal B})\otimes^{(1)} {\cal C}](((A',B'),C'),((A'',B''),C''))\otimes[({\cal A}\otimes^{(1)} {\cal B})\otimes^{(1)} {\cal C}](((A,B),C),((A',B'),C'))$$
  and the last bullet represents
  $$[{\cal A}\otimes^{(1)} ({\cal B}\otimes^{(1)} {\cal C})]((A,(B,C)),(A'',(B'',C''))).$$
  $$
  \xymatrix{
  \bullet
  \ar[rr]^{M}
  \ar[d]^{\alpha^{(1)} \otimes \alpha^{(1)}}
  &&\bullet
  \ar[d]^{\alpha^{(1)}}
  \\
  \bullet
  \ar[rr]_{M}
  &&\bullet
  }
  $$
  
  In ${\cal V}$ let $X = {\cal A}(A',A'')$, $Y = {\cal B}(B',B'')$, $Z = {\cal C}(C',C'')$, $X' = {\cal A}(A,A')$, $Y' = {\cal B}(B,B')$ and $Z' = {\cal C}(C,C')$
  Then expanding the above diagram 
  (where I leave out some parentheses for clarity and denote various composites
  of $\alpha$ and $c$ by unlabeled arrows) we have
  
  \noindent
                	          \begin{center}
    	          \resizebox{6.5in}{!}{
  $$
  \xymatrix{
  &(X\otimes Y)\otimes Z \otimes (X'\otimes Y')\otimes Z' 
  \ar[dr]
  \ar[dl]\\
  X\otimes (Y\otimes Z) \otimes X'\otimes (Y'\otimes Z')
  \ar[d]
  &&(X\otimes Y)\otimes (X'\otimes Y') \otimes Z \otimes Z'
  \ar[d]\\
  X\otimes X' \otimes (Y\otimes Z)\otimes (Y'\otimes Z')
  \ar[d]
  &&[X\otimes Y\otimes X'\otimes Y'] \otimes (Z\otimes Z')
  \ar[d]
  \\
  (X\otimes X') \otimes [Y\otimes Z\otimes Y'\otimes Z']
  \ar[d]
  &&[(X\otimes X')\otimes (Y\otimes Y')] \otimes (Z\otimes Z')
  \ar[dll]^{\alpha}
  \ar[d]^{(M\otimes M) \otimes M}\\
  (X\otimes X') \otimes [(Y\otimes Y')\otimes (Z\otimes Z')]
  \ar[dr]^{M\otimes (M \otimes M)}
  &&[{\cal A}(A,A'')\otimes {\cal B}(B,B'')]\otimes {\cal C}(C,C'')
  \ar[dl]^{\alpha}\\
  &{\cal A}(A,A'')\otimes [{\cal B}(B,B'')\otimes {\cal C}(C,C'')]
  }
  $$
  }
  \end{center}
  
  The bottom quadrilateral commutes by naturality of $\alpha$. The top region must then commute for the diagram to commute.
  These basic nodes must be present regardless of the choice of braid by which the composition morphism is defined.  
  Notice that  
  the right and left legs have the following underlying braids
  for some examples of various choices of $b.$
  The first is the one used in the original definition given above.
  
  \clearpage
   $$
   b_{(1)}=
   \xy 0;/r1pc/:
   ,{\xcapv-@(0)}+(1,1),{\vtwist}+(2,1),{\xcapv-@(0)}
   \endxy
   \text{ Functoriality follows from: }
   \xy 0;/r1pc/:
      ,{\xcapv-@(0)}+(1,1),{\xcapv-@(0)}+(1,1),{\vtwist}+(2,1),{\xcapv-@(0)}+(1,1),{\xcapv-@(0)}+(-5,0)
      ,{\xcapv-@(0)}+(1,1),{\vtwist}+(2,1),{\xcapv-@(0)}+(1,1),{\xcapv-@(0)}+(1,1),{\xcapv-@(0)}+(-5,0)
      ,{\xcapv-@(0)}+(1,1),{\xcapv-@(0)}+(1,1),{\xcapv-@(0)}+(1,1),{\vtwist}+(2,1),{\xcapv-@(0)}
   \endxy
   \xymatrix{&=}
   \xy 0;/r1pc/:
      ,{\xcapv-@(0)}+(1,1),{\xcapv-@(0)}+(1,1),{\vtwist}+(2,1),{\xcapv-@(0)}+(1,1),{\xcapv-@(0)}+(-5,0)
      ,{\xcapv-@(0)}+(1,1),{\xcapv-@(0)}+(1,1),{\xcapv-@(0)}+(1,1),{\vtwist}+(2,1),{\xcapv-@(0)}+(-5,0)
      ,{\xcapv-@(0)}+(1,1),{\vtwist}+(2,1),{\xcapv-@(0)}+(1,1),{\xcapv-@(0)}+(1,1),{\xcapv-@(0)}
   \endxy
   $$
   
   $$
   b_{(2)}=
   \xy 0;/r1pc/:
   ,{\xcapv-@(0)}+(1,1),{\xcapv-@(0)}+(1,1),{\vtwist}+(-2,0)
   ,{\xcapv-@(0)}+(1,1),{\xcapv-@(0)}+(1,1),{\vtwist}+(-2,0)
   ,{\xcapv-@(0)}+(1,1),{\vtwist}+(2,1),{\xcapv-@(0)}
   \endxy
   \text{ Functoriality follows from: }
   \xy 0;/r1pc/:
      ,{\xcapv-@(0)}+(1,1),{\xcapv-@(0)}+(1,1),{\xcapv-@(0)}+(1,1),{\vtwist}+(2,1),{\xcapv-@(0)}+(-5,0)
      ,{\xcapv-@(0)}+(1,1),{\xcapv-@(0)}+(1,1),{\xcapv-@(0)}+(1,1),{\xcapv-@(0)}+(1,1),{\vtwist}+(-4,0)
      ,{\xcapv-@(0)}+(1,1),{\xcapv-@(0)}+(1,1),{\xcapv-@(0)}+(1,1),{\xcapv-@(0)}+(1,1),{\vtwist}+(-4,0)
      ,{\xcapv-@(0)}+(1,1),{\xcapv-@(0)}+(1,1),{\xcapv-@(0)}+(1,1),{\vtwist}+(2,1),{\xcapv-@(0)}+(-5,0)
      ,{\xcapv-@(0)}+(1,1),{\xcapv-@(0)}+(1,1),{\vtwist}+(2,1),{\xcapv-@(0)}+(1,1),{\xcapv-@(0)}+(-5,0)
      ,{\xcapv-@(0)}+(1,1),{\vtwist}+(2,1),{\xcapv-@(0)}+(1,1),{\xcapv-@(0)}+(1,1),{\xcapv-@(0)}+(-5,0)
      ,{\xcapv-@(0)}+(1,1),{\xcapv-@(0)}+(1,1),{\xcapv-@(0)}+(1,1),{\xcapv-@(0)}+(1,1),{\vtwist}+(-4,0)
      ,{\xcapv-@(0)}+(1,1),{\xcapv-@(0)}+(1,1),{\xcapv-@(0)}+(1,1),{\xcapv-@(0)}+(1,1),{\vtwist}+(-4,0)
      ,{\xcapv-@(0)}+(1,1),{\xcapv-@(0)}+(1,1),{\xcapv-@(0)}+(1,1),{\vtwist}+(2,1),{\xcapv-@(0)}+(-5,0)
   \endxy
   \xymatrix{&=}
   \xy 0;/r1pc/:
      ,{\xcapv-@(0)}+(1,1),{\xcapv-@(0)}+(1,1),{\xcapv-@(0)}+(1,1),{\xcapv-@(0)}+(1,1),{\vtwist}+(-4,0)
      ,{\xcapv-@(0)}+(1,1),{\xcapv-@(0)}+(1,1),{\xcapv-@(0)}+(1,1),{\vtwist}+(2,1),{\xcapv-@(0)}+(-5,0)
      ,{\xcapv-@(0)}+(1,1),{\xcapv-@(0)}+(1,1),{\xcapv-@(0)}+(1,1),{\vtwist}+(2,1),{\xcapv-@(0)}+(-5,0)
      ,{\xcapv-@(0)}+(1,1),{\xcapv-@(0)}+(1,1),{\vtwist}+(2,1),{\vtwist}+(-4,0)
      ,{\xcapv-@(0)}+(1,1),{\xcapv-@(0)}+(1,1),{\xcapv-@(0)}+(1,1),{\vtwist}+(2,1),{\xcapv-@(0)}+(-5,0)
      ,{\xcapv-@(0)}+(1,1),{\xcapv-@(0)}+(1,1),{\vtwist}+(2,1),{\xcapv-@(0)}+(1,1),{\xcapv-@(0)}+(-5,0)
      ,{\xcapv-@(0)}+(1,1),{\xcapv-@(0)}+(1,1),{\vtwist}+(2,1),{\xcapv-@(0)}+(1,1),{\xcapv-@(0)}+(-5,0)
      ,{\xcapv-@(0)}+(1,1),{\vtwist}+(2,1),{\xcapv-@(0)}+(1,1),{\xcapv-@(0)}+(1,1),{\xcapv-@(0)}
   \endxy
   $$
  \clearpage
   $$
   b_{(3)}=
   \xy 0;/r1pc/:
   ,{\xcapv-@(0)}+(1,1),{\vtwist}+(2,1),{\xcapv-@(0)}+(-3,0)
   ,{\vtwist}+(2,1),{\xcapv-@(0)}+(1,1),{\xcapv-@(0)}+(-3,0)
   ,{\vtwist}+(2,1),{\xcapv-@(0)}+(1,1),{\xcapv-@(0)}
   \endxy
   \text{ Functoriality does not follow since: }
   \xy 0;/r1pc/:
      ,{\xcapv-@(0)}+(1,1),{\xcapv-@(0)}+(1,1),{\vtwist}+(2,1),{\xcapv-@(0)}+(1,1),{\xcapv-@(0)}+(-5,0)
      ,{\xcapv-@(0)}+(1,1),{\vtwist}+(2,1),{\xcapv-@(0)}+(1,1),{\xcapv-@(0)}+(1,1),{\xcapv-@(0)}+(-5,0)
      ,{\vtwist}+(2,1),{\xcapv-@(0)}+(1,1),{\xcapv-@(0)}+(1,1),{\xcapv-@(0)}+(1,1),{\xcapv-@(0)}+(-5,0)
      ,{\vtwist}+(2,1),{\xcapv-@(0)}+(1,1),{\xcapv-@(0)}+(1,1),{\xcapv-@(0)}+(1,1),{\xcapv-@(0)}+(-5,0)
      ,{\xcapv-@(0)}+(1,1),{\xcapv-@(0)}+(1,1),{\xcapv-@(0)}+(1,1),{\vtwist}+(2,1),{\xcapv-@(0)}+(-5,0)
      ,{\xcapv-@(0)}+(1,1),{\xcapv-@(0)}+(1,1),{\vtwist}+(2,1),{\xcapv-@(0)}+(1,1),{\xcapv-@(0)}+(-5,0)
      ,{\xcapv-@(0)}+(1,1),{\xcapv-@(0)}+(1,1),{\vtwist}+(2,1),{\xcapv-@(0)}+(1,1),{\xcapv-@(0)}+(-5,0)
   \endxy
   \xymatrix{&\ne}
   \xy 0;/r1pc/:
      ,{\xcapv-@(0)}+(1,1),{\xcapv-@(0)}+(1,1),{\vtwist}+(2,1),{\xcapv-@(0)}+(1,1),{\xcapv-@(0)}+(-5,0)
      ,{\xcapv-@(0)}+(1,1),{\vtwist}+(2,1),{\vtwist}+(2,1),{\xcapv-@(0)}+(-5,0)
      ,{\vtwist}+(2,1),{\vtwist}+(2,1),{\xcapv-@(0)}+(1,1),{\xcapv-@(0)}+(-5,0)
      ,{\xcapv-@(0)}+(1,1),{\vtwist}+(2,1),{\xcapv-@(0)}+(1,1),{\xcapv-@(0)}+(1,1),{\xcapv-@(0)}+(-5,0)
      ,{\xcapv-@(0)}+(1,1),{\vtwist}+(2,1),{\xcapv-@(0)}+(1,1),{\xcapv-@(0)}+(1,1),{\xcapv-@(0)}+(-5,0)
      ,{\vtwist}+(2,1),{\vtwist}+(2,1),{\xcapv-@(0)}+(1,1),{\xcapv-@(0)}+(-5,0)
      ,{\xcapv-@(0)}+(1,1),{\vtwist}+(2,1),{\xcapv-@(0)}+(1,1),{\xcapv-@(0)}+(1,1),{\xcapv-@(0)}+(-5,0)
      ,{\xcapv-@(0)}+(1,1),{\vtwist}+(2,1),{\xcapv-@(0)}+(1,1),{\xcapv-@(0)}+(1,1),{\xcapv-@(0)}+(-5,0)
      ,{\vtwist}+(2,1),{\xcapv-@(0)}+(1,1),{\xcapv-@(0)}+(1,1),{\xcapv-@(0)}+(1,1),{\xcapv-@(0)}+(-5,0)
      ,{\vtwist}+(2,1),{\xcapv-@(0)}+(1,1),{\xcapv-@(0)}+(1,1),{\xcapv-@(0)}+(1,1),{\xcapv-@(0)}+(-5,0)
   \endxy
   $$
   
   $$
   b_{(4)}=
   \xy 0;/r1pc/:
   ,{\xcapv-@(0)}+(1,1),{\vtwist}+(2,1),{\xcapv-@(0)}+(-3,0)
   ,{\xcapv-@(0)}+(1,1),{\vtwist}+(2,1),{\xcapv-@(0)}+(-3,0)
   ,{\xcapv-@(0)}+(1,1),{\vtwist}+(2,1),{\xcapv-@(0)}+(-3,0)
   \endxy
   \text{ Functoriality does not follow since: }
   \xy 0;/r1pc/:
      ,{\xcapv-@(0)}+(1,1),{\xcapv-@(0)}+(1,1),{\vtwist}+(2,1),{\xcapv-@(0)}+(1,1),{\xcapv-@(0)}+(-5,0)
      ,{\xcapv-@(0)}+(1,1),{\vtwist}+(2,1),{\xcapv-@(0)}+(1,1),{\xcapv-@(0)}+(1,1),{\xcapv-@(0)}+(-5,0)
      ,{\xcapv-@(0)}+(1,1),{\vtwist}+(2,1),{\xcapv-@(0)}+(1,1),{\xcapv-@(0)}+(1,1),{\xcapv-@(0)}+(-5,0)
      ,{\xcapv-@(0)}+(1,1),{\xcapv-@(0)}+(1,1),{\vtwist}+(2,1),{\xcapv-@(0)}+(1,1),{\xcapv-@(0)}+(-5,0)
      ,{\xcapv-@(0)}+(1,1),{\xcapv-@(0)}+(1,1),{\vtwist}+(2,1),{\xcapv-@(0)}+(1,1),{\xcapv-@(0)}+(-5,0)
      ,{\xcapv-@(0)}+(1,1),{\vtwist}+(2,1),{\xcapv-@(0)}+(1,1),{\xcapv-@(0)}+(1,1),{\xcapv-@(0)}+(-5,0)
      ,{\xcapv-@(0)}+(1,1),{\xcapv-@(0)}+(1,1),{\xcapv-@(0)}+(1,1),{\vtwist}+(2,1),{\xcapv-@(0)}+(-5,0)
      ,{\xcapv-@(0)}+(1,1),{\xcapv-@(0)}+(1,1),{\xcapv-@(0)}+(1,1),{\vtwist}+(2,1),{\xcapv-@(0)}+(-5,0)
      ,{\xcapv-@(0)}+(1,1),{\xcapv-@(0)}+(1,1),{\xcapv-@(0)}+(1,1),{\vtwist}+(2,1),{\xcapv-@(0)}+(-5,0)
   \endxy
   \xymatrix{&\ne}
   \xy 0;/r1pc/:
      ,{\xcapv-@(0)}+(1,1),{\xcapv-@(0)}+(1,1),{\vtwist}+(2,1),{\xcapv-@(0)}+(1,1),{\xcapv-@(0)}+(-5,0)
      ,{\xcapv-@(0)}+(1,1),{\xcapv-@(0)}+(1,1),{\xcapv-@(0)}+(1,1),{\vtwist}+(2,1),{\xcapv-@(0)}+(-5,0)
      ,{\xcapv-@(0)}+(1,1),{\xcapv-@(0)}+(1,1),{\xcapv-@(0)}+(1,1),{\vtwist}+(2,1),{\xcapv-@(0)}+(-5,0)
      ,{\xcapv-@(0)}+(1,1),{\xcapv-@(0)}+(1,1),{\vtwist}+(2,1),{\xcapv-@(0)}+(1,1),{\xcapv-@(0)}+(-5,0)
      ,{\xcapv-@(0)}+(1,1),{\xcapv-@(0)}+(1,1),{\vtwist}+(2,1),{\xcapv-@(0)}+(1,1),{\xcapv-@(0)}+(-5,0)
      ,{\xcapv-@(0)}+(1,1),{\xcapv-@(0)}+(1,1),{\xcapv-@(0)}+(1,1),{\vtwist}+(2,1),{\xcapv-@(0)}+(-5,0)
      ,{\xcapv-@(0)}+(1,1),{\vtwist}+(2,1),{\xcapv-@(0)}+(1,1),{\xcapv-@(0)}+(1,1),{\xcapv-@(0)}+(-5,0)
      ,{\xcapv-@(0)}+(1,1),{\vtwist}+(2,1),{\xcapv-@(0)}+(1,1),{\xcapv-@(0)}+(1,1),{\xcapv-@(0)}+(-5,0)
      ,{\xcapv-@(0)}+(1,1),{\vtwist}+(2,1),{\xcapv-@(0)}+(1,1),{\xcapv-@(0)}+(1,1),{\xcapv-@(0)}+(-5,0)
   \endxy
   $$
 \clearpage
  
  A comparison with the previous  examples is of interest. Braids (2) and (3) are 180 degree rotations of each other.
  Notice that the second braid in the set of functoriality examples
  leads to an equality that is actually the same as for the third 
  braid in the set of associativity examples. To see this the
  page must be rotated by 180 degrees. Similarly, the inequality preventing 
  braid (2) from being associative is the 180 degree 
  rotation of the inequality preventing braid (3) from being functorial. 
  Braid (1) is its own 180 degree rotation, and the two 
  braids proving it to be the underlying braid of an associative composition 
  morphism are the same two that show it to underlie a functorial associator .
  Braid (4) is its own 180 degree rotation, and the two 
  braids preventing it from being associative are the same two that obstruct it from being functorial. 
  Thus there is a certain kind of duality 
  between the requirements of associativity of the enriched composition and the 
  functoriality of the associator. My conjecture is that
  there are no other braids underlying the composition of a product of enriched
  categories besides the braid (1) above (and its inverse) that fulfill both obligations.
  If we were considering a strictly associative monoidal category ${\cal V}$ then 
  the condition of a functorial associator
  would become a condition of a well defined composition morphism. I think that including 
  the coherent associator is more
  enlightening. This area certainly merits more scrutiny. It may be that as the categorical delooping is more
  completely understood, information may flow the other way and we will learn new things 
  about the braid group, such as 
  an answer 
  to the conjecture. For now though these observations serve to highlight how the $k$--fold monoidal case is more
  suited to delooping than the braided case.
  
  Notice that in the symmetric case the axioms of enriched categories for ${\cal A}\otimes^{(1)} {\cal B}$ and 
  the existence of a coherent 2-natural transformation follow from the coherence 
  of symmetric categories and the 
  enriched axioms for ${\cal A}$ and ${\cal B}.$ 
  
  The unit ${\cal V}$--category ${\cal I}$ has only one object $0$ and ${\cal I}(0,0)=I$ the 
  unit in ${\cal V}$. Thus we have that, using the multiplication defined with braid (1), ${\cal V}$--Cat 
  is a monoidal 2--category. It remains to consider just why it is 
  that ${\cal V}$--Cat is braided if and only if ${\cal V}$ is symmetric, and if so, then ${\cal V}$--Cat 
  is symmetric as well.
  A braiding $c^{(1)}$ on ${\cal V}$--Cat would be a 2--natural transformation 
  $c^{(1)}_{{\cal A}{\cal B}}$ is a ${\cal V}$--functor
  ${\cal A}\otimes {\cal B} \to {\cal B}\otimes {\cal A}$. Of course $c^{(1)}_{{\cal A}{\cal B}}((A,B)) = (B,A).$
  Now to be precise we define $c^{(1)}$ to be based upon $c$ to mean that
  $$c^{(1)}_{{\cal A}{\cal B}_{(A,B)(A',B')}}: ({\cal A}\otimes^{(1)} {\cal B})((A,B),(A',B')) \to ({\cal B}\otimes^{(1)} {\cal A})((B,A),(B',A'))$$
  is exactly equal to
  $$ c_{{\cal A}(A,A'){\cal B}(B,B')}:{\cal A}(A,A')\otimes {\cal B}(B,B') \to {\cal B}(B,B')\otimes {\cal A}(A,A')$$
  This must be checked for ${\cal V}$--functoriality. 
  \clearpage
  Again the unit axioms are trivial and we consider the 
  more interesting 
  associativity of hom--object morphisms property. The following diagram must commute
  $$
  \xymatrix{
  ({\cal A}\otimes^{(1)} {\cal B})((A',B'),(A'',B''))\otimes({\cal A}\otimes^{(1)} {\cal B})((A,B),(A',B'))
  \ar[rr]^-{M}
  \ar[d]^{c^{(1)}\otimes c^{(1)}}
  &&({\cal A}\otimes^{(1)} {\cal B})((A,B),(A'',B''))
  \ar[d]^{c^{(1)}}
  \\
  ({\cal B}\otimes^{(1)} {\cal A})((B',A'),(B'',A''))\otimes ({\cal B}\otimes^{(1)} {\cal A})((B,A),(B',A'))
  \ar[rr]^-M
  &&({\cal B}\otimes^{(1)} {\cal A})((B,A),(B'',A''))
  }
  $$
  Let $X = {\cal A}(A',A'')$, $Y = {\cal B}(B',B'')$, $Z = {\cal A}(A,A')$ and $W = {\cal B}(B,B')$
  Then expanding the above diagram using the composition defined as above (denoting various composites
  of $\alpha$ by unlabeled arrows) we have
  $$
  \xymatrix{
  &(X\otimes Y) \otimes (Z\otimes W)
  \ar[dr]
  \ar[dl]|{c_{XY} \otimes  c_{ZW}}\\
  (Y\otimes X) \otimes (W\otimes Z)
  \ar[d]
  &&X\otimes ((Y \otimes Z)\otimes W)
  \ar[d]^{1 \otimes (c_{YZ} \otimes 1)}\\
  Y\otimes ((X \otimes W)\otimes Z)
  \ar[d]^{1 \otimes (c_{XW} \otimes 1)}
  &&X\otimes ((Z \otimes Y)\otimes W)
  \ar[d]
  \\
  Y\otimes ((W \otimes X)\otimes Z)
  \ar[d]
  &&(X\otimes Z) \otimes (Y\otimes W)
  \ar[dll]|{c_{(X\otimes Z)(Y\otimes W)}}
  \ar[d]^{M_{AA'A''}\otimes M_{BB'B''}}\\
  (Y\otimes W) \otimes (X\otimes Z)
  \ar[dr]|{M_{BB'B''}\otimes M_{AA'A''}}
  &&{\cal A}(A,A'')\otimes {\cal B}(B,B'')
  \ar[dl]^{c}\\
  &{\cal B}(B,B'')\otimes {\cal A}(A,A'')
  }
  $$
  The bottom quadrilateral commutes by naturality of $c$. The top region must then commute for the diagram to commute, but 
  the left and right legs have the following underlying braids
  $$
  
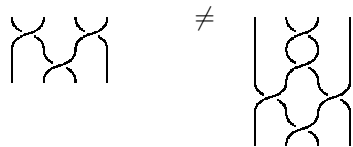
\xy 0;/r1pc/:
  ,{\vtwist}+(2,1),{\vtwist}+(-2,0)
  ,{\xcapv-@(0)}+(1,1),{\vtwist}+(2,1),{\xcapv-@(0)}+(-3,0)
  \endxy
  \xymatrix{&\ne}
  \xy 0;/r1pc/:
  ,{\xcapv-@(0)}+(1,1),{\vtwist}+(2,1),{\xcapv-@(0)}+(-3,0)
  ,{\xcapv-@(0)}+(1,1),{\vtwist}+(2,1),{\xcapv-@(0)}+(-3,0)
  ,{\vtwist}+(2,1),{\vtwist}+(-2,0)
  ,{\xcapv-@(0)}+(1,1),{\vtwist}+(2,1),{\xcapv-@(0)}+(-3,0)
  \endxy
  $$
  
  Thus neither braid (1) nor its inverse can give a monoidal structure with a braiding based on the original braiding,
  In fact, it is easy to show that composition morphisms for product enriched categories with any underlying
  braid $x$ will fail to produce a braiding in ${\cal V}$--Cat. Notice that in the above braid 
  inequality each side of the 
  inequality consists of the braid that underlies the definition of the composition morphism, 
  in this case $b_{(1)}$, and
  an additional braid that underlies the segment of the preceding diagram that corresponds to a composite of $c^{(1)}.$
  In terms of braid generators the left side of the braid inequality begins with $\sigma_1\sigma_3$ 
  corresponding to $c_{XY} \otimes  c_{ZW}$ and the right side of the braid inequality 
  ends with $\sigma_2\sigma_1\sigma_3\sigma_2$ corresponding to $c_{(X\otimes Z)(Y\otimes W)}.$ Since the same braid 
  $x$ must
  end the left side as begins the right side, then for the diagram to commute we require $x\sigma_1\sigma_3 = \sigma_2\sigma_1\sigma_3\sigma_2x.$
  This implies $\sigma_1\sigma_3 = x^{-1}\sigma_2\sigma_1\sigma_3\sigma_2x$, but for this equality to hold, $x$ must contain as a factor $\sigma_2^{-1}.$
  We see that for every factor of $\sigma_2^{-1}$ in $x$, there is of course a factor of $\sigma_2$ in $x^{-1}$, and so the
  equality $\sigma_1\sigma_3 = x^{-1}\sigma_2\sigma_1\sigma_3\sigma_2x$ can never hold due to the inability of the right side to reduce away its 
  factors of $\sigma_2$ in route to becoming the left side.
  
    It is also interesting to note that the braid inequality above is the 180 degree rotation of the one which
  implies that $({\cal A}\otimes^{(1)} {\cal B})^{op} \ne {\cal A}^{op}\otimes^{(1)} {\cal B}^{op}.$
  Thus the proof also implies that the latter inequality holds for product enriched categories with any braid $x$ 
  underlying their composition morphisms.

  It is quickly seen that if $c$ is a symmetry then in the second half of the braid inequality the upper portion
    of the braid consists of $c_{YZ}$ and $c_{ZY} = c_{YZ}^{-1}$ so in fact equality holds. In that case then the 
  derived braiding $c^{(1)}$ is a symmetry simply due to the definition.

  
  \clearpage
  \newpage
  \MySection{$k$-fold Monoidal Categories}
    
    In this section I closely follow the authors of \cite{Balt} in defining a notion of iterated monoidal category.
    For those readers familiar with that source, note that I vary from their definition only 
    by including associativity 
    up to natural coherent isomorphisms.
    This includes changing the basic picture from monoids to something that is a monoid only
    up to a monoidal natural transformation.
    We (and in this section ``we'' is not merely imperial, since so much is directly from \cite{Balt})
    start by
    defining a slightly nonstandard variant
    of monoidal functor. It is usually required in a definition of monoidal functor
    that $\eta$ be an isomorphism. The authors of \cite {Balt} note that it is crucial not to make this
    requirement. 
    
    \begin{definition}
    A {\it monoidal functor} $(F,\eta) :{\cal C}\to{\cal D}$ between monoidal categories consists
    of a functor $F$ such that $F(I)=I$ together with a natural transformation
    $$
    \eta_{AB}:F(A)\otimes F(B)\to F(A\otimes B),
    $$
    which satisfies the following conditions
    \begin{enumerate}
    \item Internal Associativity: The following diagram commutes 
    $$
    \diagram
    (F(A)\otimes F(B))\otimes F(C)
    \rrto^{\eta_{AB}\otimes 1_{F(C)}}
    \dto^{\alpha}
    &&F(A\otimes B)\otimes F(C)
    \dto^{\eta_{(A\otimes B)C}}\\
    F(A)\otimes (F(B)\otimes F(C))
    \ar[d]^{1_{F(A)}\otimes \eta_{BC}}
    &&F((A\otimes B)\otimes C)
    \ar[d]^{F\alpha}\\
    F(A)\otimes F(B\otimes C)
    \rrto^{\eta_{A(B\otimes C)}}
    &&F(A\otimes (B\otimes C))
    \enddiagram
    $$
    \item Internal Unit Conditions: $\eta_{AI}=\eta_{IA}=1_{F(A)}.$
    \end{enumerate}
    \end{definition}
    Given two monoidal functors $(F,\eta) :{\cal C}\to{\cal D}$ and $(G,\zeta) 
    :{\cal D}\to{\cal E}$,
    we define their composite to be the monoidal functor $(GF,\xi) :
    {\cal C}\to{\cal E}$, where
    $\xi$ denotes the composite
    $$
    \diagram
    GF(A)\otimes GF(B)\rrto^{\zeta_{F(A)F(B)}} 
    && G\bigl(F(A)\otimes F(B)\bigr)\rrto^{G(\eta_{AB})}
    &&GF(A\otimes B).
    \enddiagram
    $$
    It is easy to verify that $\xi$ satisfies the internal associativity condition above by subdividing the 
    necessary commuting diagram into two regions that commute by the axioms for $\eta$ and $\zeta$ respectively 
    and two that commute due to their naturality.
     $\mathbf{MonCat}$ is the monoidal category of monoidal categories and monoidal
    functors, with the usual Cartesian product as in ${\mathbf{Cat}}$. 
    
    A {\it monoidal natural transformation} $\theta:(F, \eta) \to (G, \zeta):{\cal D}\to{\cal E}$  is a 
    natural transformation $\theta: F\to G$ between the underlying ordinary functors of $F$ and $G$ such that the
    following diagram commutes
    $$
    \xymatrix{
    F(A)\otimes F(B)
    \ar[r]^{\eta}
    \ar[d]^{\theta_A \otimes \theta_B}
    &F(A\otimes B)
    \ar[d]^{\theta_{A \otimes B}}
    \\
    G(A)\otimes G(B)
    \ar[r]^{\zeta}
    &G(A\otimes B)
    }
    $$

    \clearpage

    
    \begin{definition} For our purposes a $2${\it -fold monoidal category} is a tensor object 
    in $\mathbf{MonCat}$.
    This means that we are given a monoidal category $({\cal V},\otimes_1,\alpha^1,I)$ and a 
    monoidal functor
    $(\otimes_2,\eta):{\cal V}\times{\cal V}\to{\cal V}$ which satisfies
    \begin{enumerate}
    \item External Associativity: the following diagram describes a monoidal natural transformation
    $\alpha^2$ in $\mathbf{MonCat}.$ 
    
    $$
    \xymatrix{
    {\cal V}\times{\cal V}\times{\cal V}
    \rrto^{(\otimes_2,\eta)\times 1_{\cal V}}
    \dto_(0.4){1_{\cal V}\times(\otimes_2,\eta)}
    && {\cal V}\times{\cal V}
    \dto^{(\otimes_2,\eta)}
    \ar@{=>}[dll]^{\alpha^2}
    \\
    {\cal V}\times{\cal V} 
    \rrto_{(\otimes_2,\eta)}
    &&
    {\cal V}
    }
    $$
    
    \item External Unit Conditions: the following diagram commutes in 
    $\mathbf{MonCat}$ 
    
    $$
    \diagram
    {\cal V}\times I 
    \rto^{\subseteq}
    \ddrto^{\cong}
    & {\cal V}\times{\cal V}
    \ddto^{(\otimes_2,\eta)}
    & I\times{\cal V}
    \lto_{\supseteq}
    \ddlto^{\cong}\\\\
    &{\cal V} 
    \enddiagram
    $$
    \item Coherence: The underlying natural transformation $\alpha^2$ satisfies
    the usual coherence pentagon.
       
    \end{enumerate}
    \end{definition}
    
    Explicitly this means that we are given a second associative binary operation
    $\otimes_2:{\cal V}\times{\cal V}\to{\cal V}$, for which $I$ is also a two-sided unit.
   We are also given a natural transformation
    $$
    \eta_{ABCD}: (A\otimes_2 B)\otimes_1 (C\otimes_2 D)\to
    (A\otimes_1 C)\otimes_2(B\otimes_1 D).
    $$
    The internal unit conditions give $\eta_{ABII}=\eta_{IIAB}=1_{A\otimes_2 B}$,
    while the external unit conditions give $\eta_{AIBI}=\eta_{IAIB}=1_{A\otimes_1 B}$.
    The internal associativity condition gives the commutative diagram
    
    $$
    \diagram
    ((U\otimes_2 V)\otimes_1 (W\otimes_2 X))\otimes_1 (Y\otimes_2 Z)
    \xto[rrr]^{\eta_{UVWX}\otimes_1 1_{Y\otimes_2 Z}}
    \ar[d]^{\alpha^1}
    &&&\bigl((U\otimes_1 W)\otimes_2(V\otimes_1 X)\bigr)\otimes_1 (Y\otimes_2 Z)
    \dto^{\eta_{(U\otimes_1 W)(V\otimes_1 X)YZ}}\\
    (U\otimes_2 V)\otimes_1 ((W\otimes_2 X)\otimes_1 (Y\otimes_2 Z))
    \dto^{1_{U\otimes_2 V}\otimes_1 \eta_{WXYZ}}
    &&&((U\otimes_1 W)\otimes_1 Y)\otimes_2((V\otimes_1 X)\otimes_1 Z)
    \ar[d]^{\alpha^1 \otimes_2 \alpha^1}
    \\
    (U\otimes_2 V)\otimes_1 \bigl((W\otimes_1 Y)\otimes_2(X\otimes_1 Z)\bigr)
    \xto[rrr]^{\eta_{UV(W\otimes_1 Y)(X\otimes_1 Z)}}
    &&& (U\otimes_1 (W\otimes_1 Y))\otimes_2(V\otimes_1 (X\otimes_1 Z))
    \enddiagram
    $$
    The external associativity condition gives the commutative diagram
    $$
    \diagram
    ((U\otimes_2 V)\otimes_2 W)\otimes_1 ((X\otimes_2 Y)\otimes_2 Z)
    \xto[rrr]^{\eta_{(U\otimes_2 V)W(X\otimes_2 Y)Z}}
    \ar[d]^{\alpha^2 \otimes_1 \alpha^2}
    &&& \bigl((U\otimes_2 V)\otimes_1 (X\otimes_2 Y)\bigr)\otimes_2(W\otimes_1 Z)
    \dto^{\eta_{UVXY}\otimes_2 1_{W\otimes_1 Z}}\\
    (U\otimes_2 (V\otimes_2 W))\otimes_1 (X\otimes_2 (Y\otimes_2 Z))
    \dto^{\eta_{U(V\otimes_2 W)X(Y\otimes_2 Z)}}
    &&&((U\otimes_1 X)\otimes_2(V\otimes_1 Y))\otimes_2(W\otimes_1 Z)
    \ar[d]^{\alpha^2}
    \\
    (U\otimes_1 X)\otimes_2\bigl((V\otimes_2 W)\otimes_1 (Y\otimes_2 Z)\bigr)
    \xto[rrr]^{1_{U\otimes_1 X}\otimes_2\eta_{VWYZ}}
    &&& (U\otimes_1 X)\otimes_2((V\otimes_1 Y)\otimes_2(W\otimes_1 Z))
    \enddiagram
    $$

    The authors of \cite{Balt} remark that we have natural transformations
    $$
    \eta_{AIIB}:A\otimes_1 B\to A\otimes_2 B\qquad\mbox{ and }\qquad
    \eta_{IABI}:A\otimes_1 B\to B\otimes_2 A.
    $$
    If they had insisted a 2-fold monoidal category be a tensor object in the category of monoidal categories
    and {\it strictly monoidal\/} functors, this would be equivalent to requiring that $\eta=1$.  In view
    of the above, they note that this would imply $A\otimes_1 B = A\otimes_2 B = B\otimes_1 A$ and similarly for morphisms.
        
     Joyal and Street \cite{JS} considered a 
    similar concept to Balteanu, Fiedorowicz, Schw${\rm \ddot a}$nzl and Vogt's idea of 2--fold monoidal category.  
    The former pair required the natural transformation $\eta_{ABCD}$ 
    to be an isomorphism and showed that the resulting category is naturally 
    equivalent to a braided monoidal category. As explained in \cite{Balt}, given such a category one 
    obtains an equivalent braided monoidal category by discarding one of the two 
    operations, say $\otimes_2$, and defining the commutativity isomorphism for the 
    remaining operation $\otimes_1$ to be the composite
    $$
    \diagram
    A\otimes_1 B\rrto^{\eta_{IABI}} 
    && B\otimes_2 A\rrto^{\eta_{BIIA}^{-1}}
    && B\otimes_1 A.
    \enddiagram
    $$
    
    \clearpage
    
    Just as in \cite{Balt} we now define a 2--fold monoidal functor
    between 2--fold monoidal categories $F:{\cal V}\to{\cal D}$. It is a functor together
    with two natural transformations:
    $$\lambda^1_{AB}:F(A)\otimes_1 F(B)\to F(A\otimes_1 B)$$
    $$\lambda^2_{AB}:F(A)\otimes_2 F(B)\to F(A\otimes_2 B)$$
    satisfying the same associativity and unit conditions as in the case of monoidal functors.
    In addition we require that the following hexagonal interchange diagram commutes:
    $$
    \diagram
    (F(A)\otimes_2 F(B))\otimes_1(F(C)\otimes_2 F(D))
    \xto[rrr]^{ {\eta_{F(A)F(B)F(C)F(D)}}}
    \dto^{ {\lambda^2_{AB}\otimes_1\lambda^2_{CD}}}
    &&&(F(A)\otimes_1 F(C))\otimes_2(F(B)\otimes_1 F(D))
    \dto^{ \lambda^1_{AC}\otimes_2\lambda^1_{BD}}\\
    F(A\otimes_2 B)\otimes_1 F(C\otimes_2 D)
    \dto^{ {\lambda^1_{(A\otimes_2 B)(C\otimes_2 D)}}}
    &&&F(A\otimes_1 C)\otimes_2 F(B\otimes_1 D)
    \dto^{ {\lambda^2_{(A\otimes_1 C)(B\otimes_1 D)}}}\\
    F((A\otimes_2 B)\otimes_1(C\otimes_2 D))
    \xto[rrr]^{ F(\eta_{ABCD})}
    &&& F((A\otimes_1 C)\otimes_2(B\otimes_1 D))
    \enddiagram
    $$
    
    We can now define the category {\textbf{2}$\mathbf{ -MonCat}$ of 2-fold monoidal categories and
    2-fold monoidal functors, and then define a 3-fold monoidal category as a tensor object in
    {\textbf{2}$\mathbf{ -MonCat}$.  From this point on, the iteration of this idea is straightforward
  and, paralleling the authors of \cite{Balt}, we arrive at the following definitions.
  
  \begin{definition} An $n${\it -fold monoidal category} is a category ${\cal V}$
          with the following structure. 
          \begin{enumerate}
          \item There are $n$ distinct multiplications
          $$\otimes_1,\otimes_2,\dots, \otimes_n:{\cal V}\times{\cal V}\to{\cal V}$$
          for each of which the associativity pentagon commutes
          
          \noindent
  		          \begin{center}
  	          \resizebox{5.5in}{!}{
          $$
          \xymatrix@C-=2pt@C=-30pt{
          &((U\otimes_i V)\otimes_i W)\otimes_i X \text{ }\text{ }
          \ar[rr]^{ \alpha^{i}_{UVW}\otimes_i 1_{X}}
          \ar[ddl]^{ \alpha^{i}_{(U\otimes_i V)WX}}
          &&\text{ }\text{ }(U\otimes_i (V\otimes_i W))\otimes_i X
          \ar[ddr]^{ \alpha^{i}_{U(V\otimes_i W)X}}&\\\\
          (U\otimes_i V)\otimes_i (W\otimes_i X)
          \ar[ddrr]^{ \alpha^{i}_{UV(W\otimes_i X)}}
          &&&&U\otimes_i ((V\otimes_i W)\otimes_i X)
          \ar[ddll]^{ 1_{U}\otimes_i \alpha^{i}_{VWX}}
          \\\\&&U\otimes_i (V\otimes_i (W\otimes_i X))&&&
          }
          $$
          }
  		                    \end{center}

          ${\cal V}$ has an object $I$ which is a strict unit
          for all the multiplications.
          \item For each pair $(i,j)$ such that $1\le i<j\le n$ there is a natural
          transformation
          $$\eta^{ij}_{ABCD}: (A\otimes_j B)\otimes_i(C\otimes_j D)\to
          (A\otimes_i C)\otimes_j(B\otimes_i D).$$
          \end{enumerate}
          These natural transformations $\eta^{ij}$ are subject to the following conditions:
          \begin{enumerate}
          \item[(a)] Internal unit condition: 
          $\eta^{ij}_{ABII}=\eta^{ij}_{IIAB}=1_{A\otimes_j B}$
          \item[(b)] External unit condition:
          $\eta^{ij}_{AIBI}=\eta^{ij}_{IAIB}=1_{A\otimes_i B}$
          \item[(c)] Internal associativity condition: The following diagram commutes
          
          $$
            \diagram
            ((U\otimes_j V)\otimes_i (W\otimes_j X))\otimes_i (Y\otimes_j Z)
            \xto[rrr]^{\eta^{ij}_{UVWX}\otimes_i 1_{Y\otimes_j Z}}
            \ar[d]^{\alpha^i}
            &&&\bigl((U\otimes_i W)\otimes_j(V\otimes_i X)\bigr)\otimes_i (Y\otimes_j Z)
            \dto^{\eta^{ij}_{(U\otimes_i W)(V\otimes_i X)YZ}}\\
            (U\otimes_j V)\otimes_i ((W\otimes_j X)\otimes_i (Y\otimes_j Z))
            \dto^{1_{U\otimes_j V}\otimes_i \eta^{ij}_{WXYZ}}
            &&&((U\otimes_i W)\otimes_i Y)\otimes_j((V\otimes_i X)\otimes_i Z)
            \ar[d]^{\alpha^i \otimes_j \alpha^i}
            \\
            (U\otimes_j V)\otimes_i \bigl((W\otimes_i Y)\otimes_j(X\otimes_i Z)\bigr)
            \xto[rrr]^{\eta^{ij}_{UV(W\otimes_i Y)(X\otimes_i Z)}}
            &&& (U\otimes_i (W\otimes_i Y))\otimes_j(V\otimes_i (X\otimes_i Z))
            \enddiagram
            $$
           \item[(d)] External associativity condition: The following diagram commutes
            $$
            \diagram
            ((U\otimes_j V)\otimes_j W)\otimes_i ((X\otimes_j Y)\otimes_j Z)
            \xto[rrr]^{\eta^{ij}_{(U\otimes_j V)W(X\otimes_j Y)Z}}
            \ar[d]^{\alpha^j \otimes_i \alpha^j}
            &&& \bigl((U\otimes_j V)\otimes_i (X\otimes_j Y)\bigr)\otimes_j(W\otimes_i Z)
            \dto^{\eta^{ij}_{UVXY}\otimes_j 1_{W\otimes_i Z}}\\
            (U\otimes_j (V\otimes_j W))\otimes_i (X\otimes_j (Y\otimes_j Z))
            \dto^{\eta^{ij}_{U(V\otimes_j W)X(Y\otimes_j Z)}}
            &&&((U\otimes_i X)\otimes_j(V\otimes_i Y))\otimes_j(W\otimes_i Z)
            \ar[d]^{\alpha^j}
            \\
            (U\otimes_i X)\otimes_j\bigl((V\otimes_j W)\otimes_i (Y\otimes_j Z)\bigr)
            \xto[rrr]^{1_{U\otimes_i X}\otimes_j\eta^{ij}_{VWYZ}}
            &&& (U\otimes_i X)\otimes_j((V\otimes_i Y)\otimes_j(W\otimes_i Z))
            \enddiagram
            $$
          
          \item[(e)] Finally it is required  for each triple $(i,j,k)$ satisfying
          $1\le i<j<k\le n$ that
          the giant hexagonal interchange diagram commutes.
         \end{enumerate}
          
          \noindent
  		          \begin{center}
  	          \resizebox{6.5in}{!}{
          $$
          \xymatrix@C=-103pt{
          &((A\otimes_k A')\otimes_j (B\otimes_k B'))\otimes_i((C\otimes_k C')\otimes_j (D\otimes_k D'))
          \ar[ddl]|{\eta^{jk}_{AA'BB'}\otimes_i \eta^{jk}_{CC'DD'}}
          \ar[ddr]|{\eta^{ij}_{(A\otimes_k A')(B\otimes_k B')(C\otimes_k C')(D\otimes_k D')}}
          \\\\
          ((A\otimes_j B)\otimes_k (A'\otimes_j B'))\otimes_i((C\otimes_j D)\otimes_k (C'\otimes_j D'))
          \ar[dd]|{\eta^{ik}_{(A\otimes_j B)(A'\otimes_j B')(C\otimes_j D)(C'\otimes_j D')}}
          &&((A\otimes_k A')\otimes_i (C\otimes_k C'))\otimes_j((B\otimes_k B')\otimes_i (D\otimes_k D'))
          \ar[dd]|{\eta^{ik}_{AA'CC'}\otimes_j \eta^{ik}_{BB'DD'}}
          \\\\
          ((A\otimes_j B)\otimes_i (C\otimes_j D))\otimes_k((A'\otimes_j B')\otimes_i (C'\otimes_j D'))
          \ar[ddr]|{\eta^{ij}_{ABCD}\otimes_k \eta^{ij}_{A'B'C'D'}}
          &&((A\otimes_i C)\otimes_k (A'\otimes_i C'))\otimes_j((B\otimes_i D)\otimes_k (B'\otimes_i D'))
          \ar[ddl]|{\eta^{jk}_{(A\otimes_i C)(A'\otimes_i C')(B\otimes_i D)(B'\otimes_i D')}}
          \\\\
          &((A\otimes_i C)\otimes_j (B\otimes_i D))\otimes_k((A'\otimes_i C')\otimes_j (B'\otimes_i D'))
          }
          $$
          }
  		                    \end{center}

        \end{definition}

  \vspace{.5cm}
  \begin{definition} An {\it $n$--fold monoidal functor}
  $(F,\lambda^1,\dots,\lambda^n):{\cal C}\to{\cal D}$ between $n$--fold monoidal categories
  consists of a functor $F$ such that $F(I)=I$ together with natural
  transformations
  $$\lambda^i_{AB}:F(A)\otimes_i F(B)\to F(A\otimes_i B)\quad i=1,2,\dots, n$$
  satisfying the same associativity and unit conditions as monoidal functors.
  In addition the following hexagonal interchange diagram commutes:
  $$
  \diagram
  (F(A)\otimes_j F(B))\otimes_i(F(C)\otimes_j F(D))
  \xto[rrr]^{ {\eta^{ij}_{F(A)F(B)F(C)F(D)}}}
  \dto^{ {\lambda^j_{AB}\otimes_i\lambda^j_{CD}}}
  &&&(F(A)\otimes_i F(C))\otimes_j(F(B)\otimes_i F(D))
  \dto^{ \lambda^i_{AC}\otimes_j\lambda^i_{BD}}\\
  F(A\otimes_j B)\otimes_i F(C\otimes_j D)
  \dto^{ {\lambda^i_{(A\otimes_j B)(C\otimes_j D)}}}
  &&&F(A\otimes_i C)\otimes_j F(B\otimes_i D)
  \dto^{ {\lambda^j_{(A\otimes_i C)(B\otimes_i D)}}}\\
  F((A\otimes_j B)\otimes_i(C\otimes_j D))
  \xto[rrr]^{ F(\eta^{ij}_{ABCD})}
  &&& F((A\otimes_i C)\otimes_j(B\otimes_i D))
  \enddiagram
  $$
  \end{definition}
  Composition of $n$-fold monoidal functors is defined 
  as for monoidal functors.

   The authors of \cite{Balt} point out that it is necessary to check that an $(n+1)$--fold monoidal category is
  the same thing as a tensor object in {\textbf{n}$\mathbf{ -MonCat}$, the category of
  $n$--fold monoidal categories and functors. Also as noticed in \cite{Balt}, the hexagonal interchange diagrams
  for the $(n+1)$--st monoidal operation regarded as an $n$--fold monoidal functor are what
  give rise to the giant hexagonal diagrams involving $\otimes_i$, $\otimes_j$ and $\otimes_{n+1}$.
  
  The authors of \cite{Balt}
  note that a symmetric monoidal category is $n$-fold monoidal for all $n$.  Just let
  $$\otimes_1=\otimes_2=\dots=\otimes_n=\otimes$$
  and define (associators added by myself)
  $$\eta^{ij}_{ABCD}=\alpha^{-1}\circ (1_A\otimes \alpha)\circ (1_A\otimes (c_{BC}\otimes 1_D))\circ (1_A\otimes \alpha^{-1})\circ \alpha$$
  for all $i<j$.

  Again as remarked in  \cite{Balt} Joyal and Street \cite{JS} require that the interchange natural transformations
  $\eta^{ij}_{ABCD}$ be isomorphisms. The latter pair observed that for $n\ge3$ the resulting sort of category
  is equivalent to a symmetric monoidal category. Thus as Balteanu, Fiedorowicz, 
  Schw${\rm \ddot a}$nzl and Vogt point out, the nerves of such categories
  have group completions which are infinite loop spaces rather than $n$--fold loop spaces.


    \clearpage
    \newpage
    \MySection{Categories Enriched over a $k$--fold Monoidal Category}
  The correct theory for enriching over a $k$--fold monoidal category ${\cal V}$
  may depend somewhat upon the point of view of the theorist. Here we are
  biased by the
  knowledge of research that reveals ${\cal V}$ to be precisely analogous to a
  $k$--fold loop space, as well 
  as by the observation that forming the category of categories enriched over
  ${\cal V}$ is something akin to delooping especially in the cases of braided and
  symmetric monoidal categories. 
  It turns out that if we let ourselves be guided by that intuition, what works quite well 
  is to simply consider
  categories enriched over the monoidal category given by ${\cal V}$ with $\otimes_{1}$.
  Of course the extra structure of ${\cal V}$ is very important --
  precisely when it comes to describing ${\cal V}$--Cat.
  We are ready to state the initial result. 
  
  \begin{theorem} \label{main:simple} For ${\cal V}$ a $k$--fold monoidal category ${\cal
  V}$--Cat is a $(k-1)$--fold
  monoidal 2-category. 
   \end{theorem} 
  
  \begin{example}\end {example} We begin by describing the $k=2$ case. ${\cal V}$ is
  2--fold monoidal with products $\otimes_{1}, \otimes_{2}.$ 
  ${\cal V}$--categories (which
  are the objects of ${\cal V}$--Cat) are defined as being enriched over $({\cal V}$,$\otimes_{1},\alpha^1,I)$.
  Here $\otimes_{1}$ plays the role of the product given by $\otimes$ in the axioms
  of section 1. We need to show that ${\cal V}$--Cat has a product.
  
  The unit object in ${\cal V}$--Cat is the enriched category ${\cal I}$ where $\left|{\cal I}\right| = \{0\}$ and 
  ${\cal I}(0,0) = I$. Of course $M_{000} = 1 = j_0.$
  The
  objects of the tensor ${\cal A} \otimes^{(1)}_{1}{\cal B}$ of two ${\cal
  V}$-categories
  ${\cal A}$ and ${\cal B}$ are simply pairs of objects, that is, elements
  in $\left|{\cal A}\right|\times\left|{\cal B}\right|$.   The hom--objects in ${\cal V}$ are given by
  $({\cal A} \otimes^{(1)}_{1}{\cal B})((A,B),(A',B')) = {\cal
  A}(A,A')\otimes_{2}{\cal B}(B,B')$. The composition morphisms that make
  ${\cal A} \otimes^{(1)}_{1}{\cal B}$ into a ${\cal V}$--category are
  immediately apparent as generalizations of the braided case. Recall that we are describing ${\cal
  A} \otimes^{(1)}_{1}{\cal B}$ as a 
  category enriched over ${\cal V}$ with product $\otimes_{1}$. Thus 
  \begin{small}
  $$M_{(A,B)(A',B')(A'',B'')} : ({\cal A} \otimes^{(1)}_{1}{\cal
  B})((A',B'),(A'',B''))\otimes_{1}({\cal A} \otimes^{(1)}_{1}{\cal
  B})((A,B),(A',B'))\to ({\cal A} \otimes^{(1)}_{1}{\cal B})((A,B),(A'',B''))$$
  \end{small}
  is given by
  $$
  \xymatrix{
  ({\cal A} \otimes^{(1)}_{1}{\cal B})((A',B'),(A'',B''))\otimes_{1}({\cal
  A} \otimes^{(1)}_{1}{\cal B})((A,B),(A',B'))
  \ar@{=}[d]\\
  ({\cal A}(A',A'')\otimes_{2}{\cal B}(B',B''))\otimes_{1}({\cal
  A}(A,A')\otimes_{2}{\cal B}(B,B'))
  \ar[d]_{\eta^{1,2}}\\
  ({\cal A}(A',A'')\otimes_{1}{\cal A}(A,A'))\otimes_{2}({\cal B}(B',B'')\otimes_{1}{\cal B}(B,B'))
  \ar[d]_{M_{AA'A''}\otimes_{2}M_{BB'B''}}\\
  ({\cal A}(A,A'')\otimes_{2}{\cal B}(B,B''))
  \ar@{=}[d]\\
  ({\cal A} \otimes^{(1)}_{1}{\cal B})((A,B),(A'',B''))
  }
  $$

  \begin{example}\end {example} Next we describe the $k=3$ case. ${\cal V}$ is
  3--fold monoidal with products $\otimes_{1}, \otimes_{2}$ and $\otimes_{3}$. ${\cal
  V}$--categories are defined as being enriched over $({\cal V}$,$\otimes_{1},\alpha^1,I).$
  Now ${\cal V}$--Cat has two products. The
  objects of both possible tensors ${\cal A} \otimes^{(1)}_{1}{\cal B}$ and ${\cal A} \otimes^{(1)}_{2}{\cal B}$ of two ${\cal
  V}$-categories
  ${\cal A}$ and ${\cal B}$ are elements
  in $\left|{\cal A}\right|\times\left|{\cal B}\right|$.   The hom--objects in ${\cal V}$ are given by
  $$({\cal A} \otimes^{(1)}_{1}{\cal B})((A,B),(A',B')) = {\cal A}(A,A')\otimes_{2}{\cal B}(B,B')$$
  just as in the previous case, and by
  $$({\cal A} \otimes^{(1)}_{2}{\cal B})((A,B),(A',B')) = {\cal A}(A,A')\otimes_{3}{\cal B}(B,B').$$
  
  The composition that makes $({\cal A} \otimes^{(1)}_{2}{\cal B})$ into a ${\cal V}$--category is analogous to that for 
  $({\cal A} \otimes^{(1)}_{1}{\cal B})$ but uses $\eta^{1,3}$ as its middle exchange morphism. 
  
  Now we need an interchange 2-natural transformation $\eta^{(1)1,2}$ for ${\cal V}$--Cat. 
  The family of morphisms $\eta^{(1)1,2}_{{\cal A}{\cal B}{\cal C}{\cal D}}$
  that make up a 2-natural transformation between the 2--functors  
  $\times^4{\cal V}$--Cat $:\to {\cal V}$--Cat in question 
  is a family of enriched functors. 
  Their action on objects
  is to send
  $$((A,B),(C,D)) \in \left|({\cal A} \otimes^{(1)}_{2}{\cal B}) \otimes^{(1)}_{1}({\cal C} \otimes^{(1)}_{2}{\cal D})\right|
   \text{ to } ((A,C),(B,D)) \in \left|({\cal A} \otimes^{(1)}_{1}{\cal C}) \otimes^{(1)}_{2}({\cal B} \otimes^{(1)}_{1}{\cal D})\right|.$$
  The correct construction of the family of hom--object morphisms in ${\cal V}$--Cat for each of these functors 
  is also clear.
  Noting that 
  $$[({\cal A} \otimes^{(1)}_{2} {\cal B}) \otimes^{(1)}_{1}({\cal C} \otimes^{(1)}_{2}{\cal D})](((A,B),(C,D)),((A',B'),(C',D')))$$
  $$=({\cal A} \otimes^{(1)}_{2}{\cal B})((A,B),(A',B'))\otimes_{2}({\cal C} \otimes^{(1)}_{2}{\cal D})((C,D),(C',D')) $$
  $$=({\cal A}(A,A')\otimes_{3}{\cal B}(B,B'))\otimes_{2}({\cal C}(C,C')\otimes_{3}{\cal D}(D,D'))$$
  and similarly 
  $$[({\cal A} \otimes^{(1)}_{1}{\cal C}) \otimes^{(1)}_{2}({\cal B} \otimes^{(1)}_{1}{\cal D})](((A,C),(B,D)),((A',C'),(B',D')))$$
  $$=({\cal A}(A,A')\otimes_{2}{\cal C}(C,C'))\otimes_{3}({\cal B}(B,B')\otimes_{2}{\cal D}(D,D'))$$
  we make the obvious identification, where by obvious I mean based upon the corresponding structure in 
  ${\cal V}$ as described earlier in the discussion of braided ${\cal V}.$ Here ``based upon'' is more freely 
  interpreted as also allowing a shift in index.
  Thus we write:
  $$\eta^{(1)1,2}_{{{\cal A}{\cal B}{\cal C}{\cal D}}_{(ABCD)(A'B'C'D')}}
  = \eta^{2,3}_{{\cal A}(A,A'){\cal B}(B,B'){\cal C}(C,C'){\cal D}(D,D')}.$$
  
  Much needs to be verified. Existence and coherence of required natural transformations, satisfaction of 
  enriched axioms and of $k$--fold monoidal axioms all must be checked. These will be dealt with next in a 
  more general setting.

  \begin{proof} of Theorem \ref{main:simple} As in the examples, ${\cal V}$--Cat is made
  up of categories enriched over $({\cal V}$,$\otimes_{1},\alpha^1,I).$ Here we define products $ \otimes^{(1)}_{1} ...  \otimes^{(1)}_{k-1}$
  in ${\cal V}$--Cat for ${\cal V}$ $k$--fold monoidal.
   We check that
  our products do make ${\cal A} \otimes^{(1)}_{2}{\cal B}$ into a ${\cal V}$--category.
  Then we check that ${\cal V}$--Cat has the 
  required coherent 2--natural transformations of
  associativity and units. We then define interchange 
  2--natural transformations $\eta^{(1)i,j}$ and check that the interchange transformations are
  2--natural and obey all the axioms required of them. It is informative to observe how these axioms are satisfied based 
  upon the axioms that ${\cal V}$ itself satisfies. It is here that we should look carefully for the algebraic reflection of 
  the topological functor $\Omega.$
  
  Again, the unit object in ${\cal V}$--Cat is the enriched category ${\cal I}$ where $\left|{\cal I}\right| = \{0\}$ and 
  ${\cal I}(0,0) = I$.
  For ${\cal V}$ $k$--fold monoidal we define the $i$th product of ${\cal
  V}$--categories ${\cal A} \otimes^{(1)}_{i}{\cal B}$
   to have objects $\in \left|{\cal A}\right|\times \left|{\cal B}\right|$
   and to have hom--objects in ${\cal V}$ given by
   
   $$({\cal A} \otimes^{(1)}_{i} {\cal B})((A,B),(A',B')) = {\cal A}(A,A')\otimes_{i+1} {\cal B}(B,B').$$
   
   Immediately we see that ${\cal V}$--Cat is $(k-1)$--fold monoidal by definition.
   The composition morphisms are 
   \begin{small}
   $$M_{(A,B)(A',B')(A'',B'')} : ({\cal A} \otimes^{(1)}_{i}{\cal
   B})((A',B'),(A'',B''))\otimes_{1}({\cal A} \otimes^{(1)}_{i}{\cal
   B})((A,B),(A',B'))\to ({\cal A} \otimes^{(1)}_{i}{\cal B})((A,B),(A'',B''))$$
   \end{small}
   given by
   $$
   \xymatrix{
   ({\cal A} \otimes^{(1)}_{i}{\cal B})((A',B'),(A'',B''))\otimes_{1}({\cal
   A} \otimes^{(1)}_{i}{\cal B})((A,B),(A',B'))
   \ar@{=}[d]\\
   ({\cal A}(A',A'')\otimes_{i+1}{\cal B}(B',B''))\otimes_{1}({\cal
   A}(A,A')\otimes_{i+1}{\cal B}(B,B'))
   \ar[d]_{\eta^{1,i+1}}\\
   ({\cal A}(A',A'')\otimes_{1}{\cal A}(A,A'))\otimes_{i+1}({\cal
   B}(B',B")\otimes_{1}{\cal B}(B,B'))
   \ar[d]_{M_{AA'A''}\otimes_{2}M_{BB'B''}}\\
   ({\cal A}(A,A'')\otimes_{i+1}{\cal B}(B,B''))
   \ar@{=}[d]\\
   ({\cal A} \otimes^{(1)}_{i}{\cal B})((A,B),(A'',B''))
  }
  $$
   The identity element is given by $j_{(A,B)} = \xymatrix{I = I \otimes_{i+1} I
   							\ar[d]^{j_A \otimes_{i+1} j_B}
   							\\{\cal A}(A,A)\otimes_{i+1} {\cal B}(B,B)
   							\ar@{=}[d]
   							\\({\cal A} \otimes^{(1)}_{i}{\cal B})((A,B),(A,B))
   							}$
   
   That each product $\otimes^{(1)}_i$ thus defined is a 2--bi--functor  ${\cal V}$--Cat $\times$ ${\cal V}$--Cat $\to$ ${\cal V}$--Cat
          is seen easily. Its action on enriched functors and natural transformations is to form formal products using $\otimes_{i+1}$ of their
          associated morphisms. That the result is a valid enriched functor or natural transformation always follows from the
       naturality of $\eta.$
   
  We first check that ${\cal A}\otimes^{(1)}_{i}{\cal B}$ is indeed properly
   enriched over ${\cal V}.$ Our definition of $M$ must obey the axioms for associativity and respect of the unit.
   For associativity the following diagram must commute, where the initial bullet represents
   $$
  [({\cal A}\otimes^{(1)}_{i}{\cal B})((A'',B''),(A''',B'''))\otimes_{1} ({\cal A}\otimes^{(1)}_{i}{\cal B})((A',B'),(A'',B''))]\otimes_{1}({\cal A}\otimes^{(1)}_{i}{\cal B})((A,B),(A',B'))
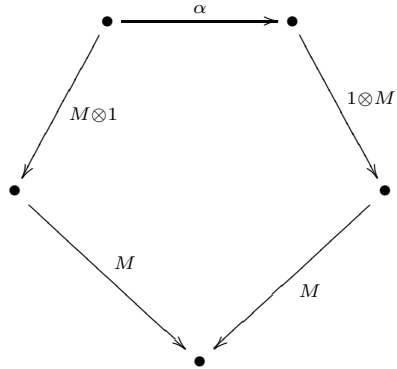
  $$

    $$
    \xymatrix{
    &\bullet
    \ar[rr]^{ \alpha}
    \ar[ddl]^{ M \otimes 1}
    &&\bullet
    \ar[ddr]^{ 1 \otimes M}&\\\\
    \bullet
    \ar[ddrr]^{ M}
    &&&&\bullet
    \ar[ddll]^{ M}
    \\\\&&\bullet
    }$$

  In the expanded diagram let
  $X={\cal A}(A,A')$, $X'={\cal A}(A',A'')$, $X''={\cal A}(A'',A''')$, $Y={\cal B}(B,B')$, $Y'={\cal B}(B',B'')$ and $Y''={\cal B}(B'',B''').$
  The exterior of the following expanded diagram is required to commute
  \clearpage
  
  \noindent
         \begin{center}
	 \resizebox{5in}{!}
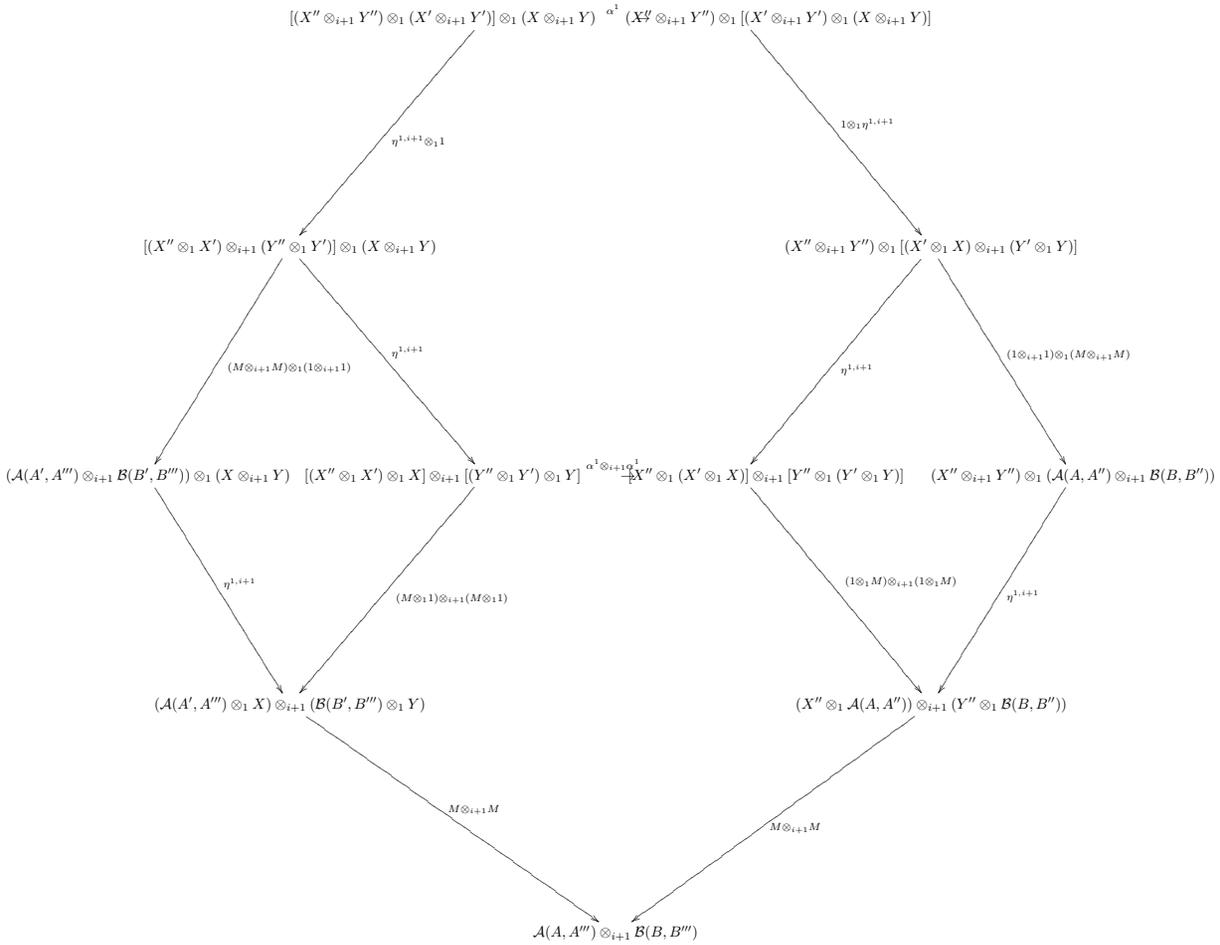
{
  \begin{sideways}
  $$
  \xymatrix@C=-105pt@R=42pt{
  &&[(X''\otimes_{i+1} Y'')\otimes_{1} (X'\otimes_{i+1} Y')]\otimes_{1} (X\otimes_{i+1} Y) \text{ }\text{ }\text{ }\text{ }\text{ }\text{ }\text{ }\text{ \normalsize{$\to$}} \text{ }\text{ }\text{ }\text{ }\text{ }\text{ }
  \ar[dddl]^{\eta^{1,i+1}\otimes_{1} 1}
  \ar@{}[rr]^{\alpha^{1}}
  &&\text{ }\text{ }\text{ }\text{ }\text{ }\text{ }\text{ }\text{ }\text{ }\text{ }\text{ }\text{ }\text{ }\text{ }\text{ }(X''\otimes_{i+1} Y'')\otimes_{1} [(X'\otimes_{i+1} Y')\otimes_{1} (X\otimes_{i+1} Y)]
  \ar[dddr]^{1 \otimes_{1} \eta^{1,i+1}}
  \\\\\\
  &[(X''\otimes_{1} X')\otimes_{i+1} (Y''\otimes_{1} Y')]\otimes_{1} (X\otimes_{i+1} Y)
  \ar[dddl]^{\text{  }(M\otimes_{i+1} M)\otimes_{1} (1\otimes_{i+1} 1)}
  \ar[dddr]^{\eta^{1,i+1}}
  &&&&(X''\otimes_{i+1} Y'')\otimes_{1} [(X'\otimes_{1} X)\otimes_{i+1} (Y'\otimes_{1} Y)]
  \ar[dddl]^{\eta^{1,i+1}}
  \ar[dddr]^{(1\otimes_{i+1} 1)\otimes_{1} (M\otimes_{i+1} M)}
  \\\\\\
  ({\cal A}(A',A''')\otimes_{i+1} {\cal B}(B',B'''))\otimes_{1} (X\otimes_{i+1} Y)
  \ar[dddr]^{\eta^{1,i+1}}
  &&[(X''\otimes_{1} X')\otimes_{1} X]\otimes_{i+1} [(Y''\otimes_{1} Y')\otimes_{1} Y] \text{ }\text{ }\text{ }\text{ }\text{ }\text{ }\text{ }\text{ \normalsize{$\to$}} \text{ }\text{ }\text{ }\text{ }\text{ }\text{ }
  \ar[dddl]^{\text{  }(M\otimes_{1} 1)\otimes_{i+1} (M\otimes_{1} 1)}
  \ar@{}[rr]^{\alpha^{1}\otimes_{i+1} \alpha^{1}}
  &&\text{ }\text{ }\text{ }\text{ }\text{ }\text{ }\text{ }\text{ }\text{ }\text{ }[X''\otimes_{1} (X'\otimes_{1} X)]\otimes_{i+1} [Y''\otimes_{1} (Y'\otimes_{1} Y)]
  \ar[dddr]^{\text{  }(1\otimes_{1} M)\otimes_{i+1} (1\otimes_{1} M)}
  &&(X''\otimes_{i+1} Y'')\otimes_{1} ({\cal A}(A,A'')\otimes_{i+1} {\cal B}(B,B''))
  \ar[dddl]^{\eta^{1,i+1}}
  \\\\\\
  &({\cal A}(A',A''')\otimes_{1} X)\otimes_{i+1} ({\cal B}(B',B''')\otimes_{1} Y)
  \ar[dddrr]^{M\otimes_{i+1} M}
  &&&&(X''\otimes_{1} {\cal A}(A,A''))\otimes_{i+1} (Y''\otimes_{1} {\cal B}(B,B''))
  \ar[dddll]^{M\otimes_{i+1} M}
  \\\\\\
  &&&{\cal A}(A,A''')\otimes_{i+1} {\cal B}(B,B''')
  } 
  $$
  \end{sideways}
  }
  \end{center}
  
  
  The lower pentagon commutes since it is two copies of the associativity axiom--one for ${\cal A}$ and one for  
  ${\cal B}.$  The two diamonds commute by the naturality of $\eta.$ The upper hexagon commutes by the internal associativity
  of $\eta.$
 \clearpage
 
 For the unit axioms we have the following compact diagram
  
  \noindent
              		          \begin{center}
	          \resizebox{6.5in}{!}{
  $$
      \xymatrix@C=-30pt{
      I\otimes_1 ({\cal A} \otimes^{(1)}_{i}{\cal B})((A,B),(A',B'))
      \ar[rrd]^{=}
      \ar[dd]_{j_{(A',B')}\otimes_1 1}
      &&&&({\cal A} \otimes^{(1)}_{i}{\cal B})((A,B),(A',B'))\otimes_1 I 
      \ar[dd]^{1\otimes_1 j_{(A,B)}}
      \ar[lld]^{=}\\
      &&({\cal A} \otimes^{(1)}_{i}{\cal B})((A,B),(A',B'))\\
      ({\cal A} \otimes^{(1)}_{i}{\cal B})((A',B'),(A',B'))\otimes_1 ({\cal A} \otimes^{(1)}_{i}{\cal B})((A,B),(A',B'))
      \ar[rru]|{\text{ }\text{ }M_{(A,B)(A',B')(A',B')}}
      &&&&({\cal A} \otimes^{(1)}_{i}{\cal B})((A,B),(A',B'))\otimes_1 ({\cal A} \otimes^{(1)}_{i}{\cal B})((A,B),(A,B))
      \ar[llu]|{\text{ }\text{ }M_{(A,B)(A,B)(A',B')}}
      }
   $$ 
  }
  \end{center}
  
  I expand the left triangle, abbreviating $X={\cal A}(A,A')$, $Y={\cal A}(A',A')$, $Z={\cal B}(B,B')$ and 
  $W={\cal B}(B',B').$ The exterior of the following must commute
  $$
  \xymatrix{
  I\otimes_1 (X\otimes_{i+1} Z)\text{ }\text{ }\text{ }
  \ar[dd]_{=}
  \ar[dr]^{=}
  \\
  &(I\otimes_1 X)\otimes_{i+1} (I\otimes_1 Z)
  \ar[dr]^{=}
  \ar[dd]_{(j_{A'} \otimes_1 1)\otimes_{i+1} (j_{B'}\otimes_1 1)}
  \\
  (I\otimes_{i+1} I)\otimes_1 (X\otimes_{i+1} Z)\text{ }\text{ }\text{ }
  \ar[dd]_{(j_{A'} \otimes_{i+1} j_{B'})\otimes_1 (1\otimes_{i+1} 1)}
  \ar[ur]^{\eta^{1,i+1}_{IIXZ}}
  &&(X\otimes_{i+1} Z)
  \\
  &(Y\otimes_1 X)\otimes_{i+1} (W\otimes_1 Z)
  \ar[ur]^{M\otimes_{i+1} M}
  \\
  (Y\otimes_{i+1} W)\otimes_1 (X\otimes_{i+1} Z)\text{ }\text{ }\text{ }
  \ar[ur]^{\eta^{1,i+1}_{YWXZ}}
  }
  $$
  
  The parallelogram commutes by naturality of $\eta$, the rightmost triangle by the 
  unit axioms of the individual ${\cal V}$--categories, and the top triangle by the internal
  unit condition for $\eta.$ The right triangle in the axiom is checked similarly.
  
  On a related note, we need to check 
  that ${\cal I} \otimes^{(1)}_i {\cal A}$ = ${\cal A}$ The object sets
  and hom--objects of the two categories in question are clearly 
  equivalent. What needs to be checked is that the composition morphisms are the same.
  Note that the composition given by
  $$
   \xymatrix{
   ({\cal I} \otimes^{(1)}_{i}{\cal A})((0,A'),(0,A''))\otimes_{1} ({\cal I} \otimes^{(1)}_{i} {\cal A})((0,A),(0,A'))
   \ar@{=}[d]\\
   (I\otimes_{i+1} {\cal A}(A',A''))\otimes_{1} (I\otimes_{i+1} {\cal A}(A,A'))
   \ar[d]_{\eta^{1,i+1}_{I{\cal A}(A',A'')I{\cal A}(A,A')}}\\
   (I\otimes_{1}I)\otimes_{i+1}({\cal A}(A',A'')\otimes_{1}{\cal A}(A,A'))
   \ar[d]_{1\otimes_{i+t} M_{AA'A''}}\\
   (I\otimes_{i+1}{\cal A}(A,A''))
   \ar@{=}[d]\\
   ({\cal I} \otimes^{(1)}_{i} {\cal A})((0,A),(0,A''))
  }
  $$
  is equivalent to simply $M_{AA'A''}$ by the external unit condition for $\eta.$
   
   Associativity in ${\cal V}$--Cat
   must hold for each $\otimes^{(1)}_{i}$. The components of 2--natural
   isomorphism 
   $$\alpha^{(1)i}_{{\cal A}{\cal B}{\cal C}}: ({\cal A} \otimes^{(1)}_{i} {\cal B})
    \otimes^{(1)}_{i} {\cal C} \to {\cal A} \otimes^{(1)}_{i} ({\cal B} \otimes^{(1)}_{i}
   {\cal C})$$
   are ${\cal V}$--functors
   that send ((A,B),C) to (A,(B,C)) and whose hom-components 
   \begin{small}
   $$\alpha^{(1)i}_{{{\cal A}{\cal B}{\cal C}}_{((A,B),C)((A',B'),C')}}: [({\cal A} \otimes^{(1)}_{i} {\cal B}) \otimes^{(1)}_{i} {\cal C}](((A,B),C),((A',B'),C'))
   \to [{\cal A} \otimes^{(1)}_{i} ({\cal B} \otimes^{(1)}_{i} {\cal C})]((A,(B,C)),(A',(B',C')))$$
   \end{small}
   are given by 
   $$\alpha^{(1)i}_{{{\cal A}{\cal B}{\cal C}}_{((A,B),C)((A',B'),C')}}
   = \alpha^{i+1}_{{\cal A}(A,A'){\cal B}(B,B'){\cal C}(C,C')}.$$
   This guarantees that the 2--natural isomorphism $\alpha^{(1)i}$ is coherent. The commutativity of the pentagon for
   the objects is trivial, and the commutativity of the pentagon for the hom--object morphisms follows directly from the 
   commutativity of the pentagon for $\alpha^{i+1}.$
  
   In order to be a functor the associator components must
   satisfy the commutativity of the
   diagrams in Definition \ref{enriched:funct}.
   \begin{enumerate}
   \item 
  
   $$
     \xymatrix{
     &\bullet
     \ar[rr]^{M}
     \ar[d]^{\alpha^{(1)i} \otimes \alpha^{(1)i}}
     &&\bullet
     \ar[d]^{\alpha^{(1)i}}&\\
     &\bullet
     \ar[rr]^{M}
     &&\bullet
     }
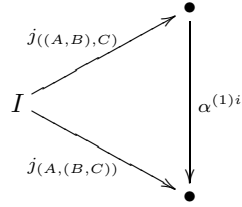
   $$
   \item
   $$
     \xymatrix{
     &&\bullet
     \ar[dd]^{ \alpha^{(1)i}}\\
     I
     \ar[rru]^{j_{((A,B),C)}}
     \ar[rrd]_{ j_{(A,(B,C))}}\\
     &&\bullet
     }
   $$
   \end{enumerate}
   
   Expanding the first using
   the definitions just given we have that the initial position in the diagram is 
   
   $$
   [({\cal A}\otimes^{(1)}_{i}{\cal B})\otimes^{(1)}_{i}{\cal C}](((A',B'),C'),((A'',B''),C''))\otimes_{1}[({\cal A}\otimes^{(1)}_{i}{\cal B})\otimes^{(1)}_{i}{\cal C}](((A,B),C),((A',B'),C'))
  $$
  
  $$
   = [({\cal A}(A',A'')\otimes_{i+1}{\cal B}(B',B''))\otimes_{i+1}{\cal C}(C',C'')]\otimes_{1}[({\cal A}(A,A')\otimes_{i+1}{\cal B}(B,B'))\otimes_{i+1}{\cal C}(C,C')]
   $$
   We let $X={\cal A}(A',A'')$, $Y={\cal B}(B',B'')$, $Z={\cal C}(C',C'')$, $X'={\cal A}(A,A')$, $Y'={\cal B}(B,B')$ and $Z'={\cal C}(C,C').$
   Then expanding the diagram we have, with an added interior arrow
   
   \noindent
               		          \begin{center}
	          \resizebox{6.5in}{!}
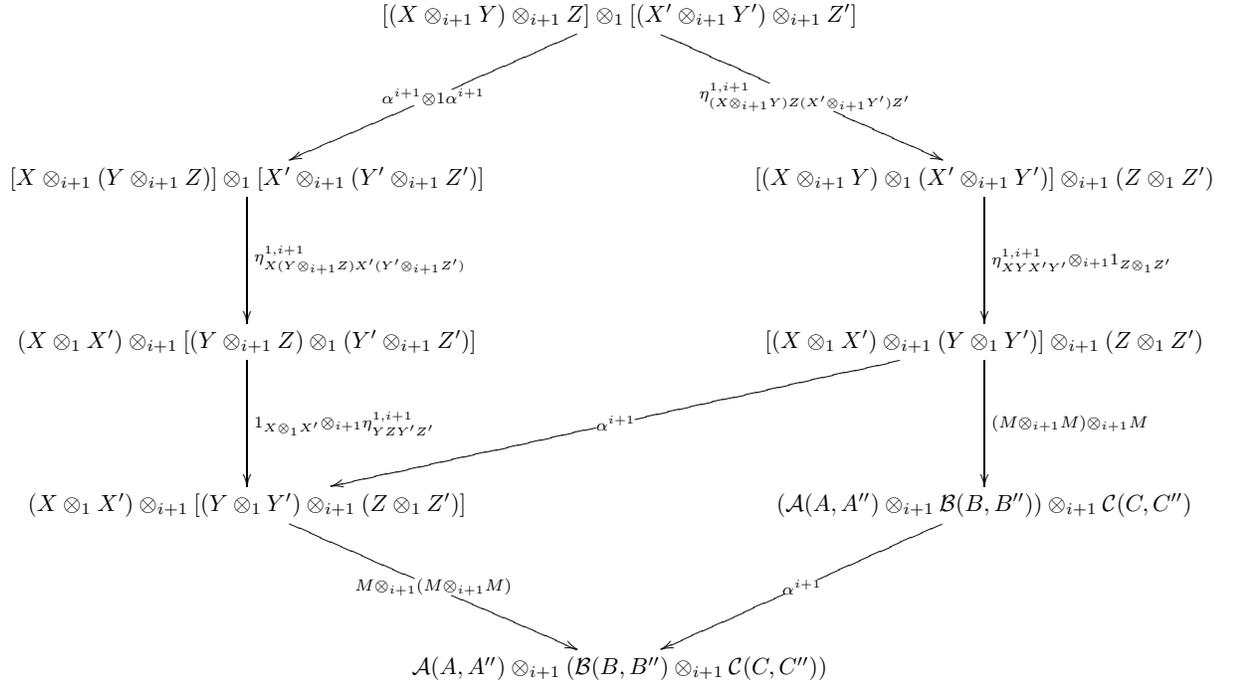
{
   $$
   \xymatrix@C=-50pt{
   &[(X\otimes_{i+1} Y)\otimes_{i+1} Z]\otimes_{1} [(X'\otimes_{i+1} Y')\otimes_{i+1} Z']
   \ar[ddr]|{\text{  }\eta^{1,i+1}_{(X\otimes_{i+1} Y)Z(X'\otimes_{i+1} Y')Z'}}
   \ar[ddl]|{\alpha^{i+1}\otimes{1} \alpha^{i+1}}\\\\
   [X\otimes_{i+1} (Y\otimes_{i+1} Z)]\otimes_{1} [X'\otimes_{i+1} (Y'\otimes_{i+1} Z')]
   \ar[dd]^{\eta^{1,i+1}_{X(Y\otimes_{i+1} Z)X'(Y'\otimes_{i+1} Z')}}
   &&[(X\otimes_{i+1} Y)\otimes_{1} (X'\otimes_{i+1} Y')]\otimes_{i+1} (Z\otimes_{1} Z')
   \ar[dd]^{\eta^{1,i+1}_{XYX'Y'}\otimes_{i+1} 1_{Z\otimes_{1} Z'}}\\\\
   (X\otimes_{1} X') \otimes_{i+1} [(Y\otimes_{i+1} Z)\otimes_{1} (Y'\otimes_{i+1} Z')]
   \ar[dd]^{1_{X\otimes_{1} X'} \otimes_{i+1} \eta^{1,i+1}_{YZY'Z'}}
   &&[(X\otimes_{1} X')\otimes_{i+1} (Y\otimes_{1} Y')]\otimes_{i+1} (Z\otimes_{1} Z')
   \ar[dd]^{(M\otimes_{i+1} M)\otimes_{i+1} M}
   \ar[ddll]|{\alpha^{i+1}}
   \\\\
   (X\otimes_{1} X')\otimes_{i+1} [(Y\otimes_{1} Y')\otimes_{i+1} (Z\otimes_{1} Z')]
   \ar[ddr]|{M\otimes_{i+1} (M\otimes_{i+1} M)}
   &&({\cal A}(A,A'')\otimes_{i+1} {\cal B}(B,B''))\otimes_{i+1}{\cal C}(C,C'')
   \ar[ddl]|{\alpha^{i+1}}\\\\
   &{\cal A}(A,A'')\otimes_{i+1} ({\cal B}(B,B'')\otimes_{i+1}{\cal C}(C,C''))
   }
   $$
   }
   \end{center}
  
   The lower quadrilateral commutes by naturality of $\alpha$, and the upper 
   hexagon commutes by the external associativity of $\eta.$
   
    The uppermost position in the expanded version of diagram number (2) is
   $$
   [({\cal A}\otimes^{(1)}_{i}{\cal B})\otimes^{(1)}_{i}{\cal C}](((A,B),C),((A,B),C))
   $$
   $$
   = [({\cal A}(A,A)\otimes_{i+1}{\cal B}(B,B))\otimes_{i+1}{\cal C}(C,C)]
   $$
   The expanded diagram is easily seen to commute by the naturality of $\alpha.$
   
   The 2--naturality of $\alpha^{(1)}$ is essentially just the naturality of its 
   components, but I think it ought to be expounded upon. Since the components of $\alpha^{(1)}$ are
   ${\cal V}$--functors the whisker diagrams for the definition of 2--naturality are defined by the 
   whiskering in ${\cal V}$--Cat. Given an arbitrary 2--cell in $\times^3{\cal V}$--Cat, i.e. 
   $(\beta,\gamma,\rho):(Q,R,S)\to (Q',R',S'):({\cal A},{\cal B},{\cal C})\to({\cal A'},{\cal B'},{\cal C'})$
   the diagrams whose composition must be equal are:
  
   $$
    \xymatrix@R-=3pt{
    &\ar@{=>}[dd]^{(\beta\otimes^{(1)}_i \gamma)\otimes^{(1)}_i \rho}\\
     ({\cal A}\otimes^{(1)}_i {\cal B})\otimes^{(1)}_i {\cal C}
    \ar@/^2pc/[rrr]^{(Q\otimes^{(1)}_i R)\otimes^{(1)}_i S}
    \ar@/_2pc/[rrr]_{(Q'\otimes^{(1)}_i R')\otimes^{(1)}_i S'}
    &&&({\cal A}'\otimes^{(1)}_i {\cal B}')\otimes^{(1)}_i {\cal C}'
    \ar[rr]^{\alpha^{(1)i}_{{\cal A}'{\cal B}'{\cal C}'}}
    &&{\cal A}'\otimes^{(1)}_i ({\cal B}'\otimes^{(1)}_i {\cal C}')
    \\
   &\\
   }
   $$
   
   $$
  = \xymatrix@R-=3pt{
   &&&\ar@{=>}[dd]^{\beta\otimes^{(1)}_i (\gamma\otimes^{(1)}_i \rho)}\\
   ({\cal A}\otimes^{(1)}_i {\cal B})\otimes^{(1)}_i {\cal C}
   \ar[rr]^{\alpha^{(1)i}_{{\cal A}{\cal B}{\cal C}}}
   &&{\cal A}\otimes^{(1)}_i ({\cal B}\otimes^{(1)}_i {\cal C})
   \ar@/^2pc/[rrr]^{Q\otimes^{(1)}_i (R\otimes^{(1)}_i S)}
   \ar@/_2pc/[rrr]_{Q'\otimes^{(1)}_i (R'\otimes^{(1)}_i S')}
   &&&{\cal A}'\otimes^{(1)}_i ({\cal B}'\otimes^{(1)}_i {\cal C}')
   \\
   &&&&\\
   }
  $$
  This is quickly seen to hold when we translate using the definitions of whiskering in ${\cal V}$--Cat, as follows.
  The $ABCD$ components of the new 2--cells are given by the exterior legs of the following diagram. They are equal 
  by naturality of $\alpha^{i+1}$ and Mac Lane's coherence theorem.
  
  \noindent
              		          \begin{center}
	          \resizebox{6.5in}{!}{
  $$
  \xymatrix{
  &I
  \ar[dr]^{=}
  \ar[dl]^{=}
  \\
  (I\otimes_{i+1} I)\otimes_{i+1} I
  \ar[rr]^{\alpha^{i+1}}
  \ar[dd]^{(\beta_A\otimes_{i+1} \gamma_B\otimes_{i+1}) \rho_C}
  &&I\otimes_{i+1} (I\otimes_{i+1} I)
  \ar[dd]^{\beta_A\otimes_{i+1} (\gamma_B\otimes_{i+1} \rho_C)}
  \\\\
  ({\cal A}'(QA,Q'A)\otimes_{i+1} {\cal B}'(RB,R'B))\otimes_{i+1} {\cal C}'(SC,S'C)
  \ar[rr]^{\alpha^{i+1}}
  &&{\cal A}'(QA,Q'A)\otimes_{i+1} ({\cal B}'(RB,R'B)\otimes_{i+1} {\cal C}'(SC,S'C))
  }
  $$
  }
  \end{center}
  
  Now we turn to consider the existence and behavior of interchange 2--natural transformations $\eta^{(1)ij}$
  for $j\ge i+1$.
  As in the example, we define the component morphisms $\eta^{(1)i,j}_{{\cal A}{\cal B}{\cal C}{\cal D}}$
  that make a 2--natural transformation between 2--functors. Each component must be an enriched functor.
  Their action on objects
  is to send $((A,B),(C,D)) \in \left|({\cal A}\otimes^{(1)}_{j} {\cal B})\otimes^{(1)}_{i} ({\cal C}\otimes^{(1)}_{j} {\cal D})\right|$
   to $((A,C),(B,D)) \in \left|({\cal A}\otimes^{(1)}_{i} {\cal C})\otimes^{(1)}_{j} ({\cal B}\otimes^{(1)}_{i} {\cal D})\right|$.
  The hom--object morphisms are given by
   $$\eta^{(1)i,j}_{{{\cal A}{\cal B}{\cal C}{\cal D}}_{(ABCD)(A'B'C'D')}} =
   \eta^{i+1,j+1}_{{\cal A}(A,A'){\cal B}(B,B'){\cal C}(C,C'){\cal D}(D,D')}.$$
  
  For this designation of $\eta^{(1)}$ to define a valid ${\cal V}$--functor, it must obey the axioms for compatibility 
  with composition and units. 
  We need commutativity of the following diagram, where the first bullet represents

  \begin{scriptsize}
  $$[({\cal A}\otimes^{(1)}_{j} {\cal B})\otimes^{(1)}_{i} ({\cal C}\otimes^{(1)}_{j} {\cal D})](((A',B'),(C',D')),((A'',B''),(C'',D'')))\otimes_1 [({\cal A}\otimes^{(1)}_{j} {\cal B})\otimes^{(1)}_{i} ({\cal C}\otimes^{(1)}_{j} {\cal D})](((A,B),(C,D)),((A',B'),(C',D')))$$
  \end{scriptsize}

  and the last bullet represents
  \begin{small}
  $$[({\cal A}\otimes^{(1)}_{i} {\cal C})\otimes^{(1)}_{j} ({\cal B}\otimes^{(1)}_{i} {\cal D})](((A,C),(B,D)),((A'',C''),(B'',D''))).$$
  \end{small}
  $$
  \xymatrix{
  \bullet
  \ar[rr]^{M}
  \ar[d]_{\eta^{(1)i,j} \otimes_1 \eta^{(1)i,j}}
  &&\bullet
  \ar[d]^{\eta^{(1)i,j}}
 \\
 \bullet
  \ar[rr]_{M}
  &&\bullet
  }
  $$
  
  If we let $X={\cal A}(A,A')$, $Y={\cal B}(B,B')$, $Z={\cal C}(C,C')$, $W={\cal D}(D,D')$, $X'={\cal A}(A',A'')$, $Y'={\cal B}(B',B'')$, $Z'={\cal C}(C',C'')$ and $W'={\cal D}(D',D'')$
  then the expanded diagram is as follows. The exterior must commute.
  \clearpage
  
  \noindent
              		          \begin{center}
	          \resizebox{4.5in}{!}{
  $$
  \begin{sideways}
  \xymatrix@C=-140pt@R=30pt{
  \text{ }
  \\
  &[(X'\otimes_{j+1} Y')\otimes_{i+1} (Z'\otimes_{j+1} W')]\otimes_{1} [(X\otimes_{j+1} Y)\otimes_{i+1} (Z\otimes_{j+1} W)]
  \ar[ddr]|{\eta^{1,i+1}_{(X'\otimes_{j+1} Y')(Z'\otimes_{j+1} W')(X\otimes_{j+1} Y)(Z\otimes_{j+1} W)}}
  \ar[ddl]|{\eta^{i+1,j+1}_{X'Y'Z'W'}\otimes_1 \eta^{i+1,j+1}_{XYZW}}\\\\
  [(X'\otimes_{i+1} Z')\otimes_{j+1} (Y'\otimes_{i+1} W')]\otimes_{1} [(X\otimes_{i+1} Z)\otimes_{j+1} (Y\otimes_{i+1} W)]
  \ar[dd]|{\eta^{1,j+1}_{(X'\otimes_{i+1} Z')(Y'\otimes_{i+1} W')(X\otimes_{i+1} Z)(Y\otimes_{i+1} W)}}
  &&[(X'\otimes_{j+1} Y')\otimes_{1} (X\otimes_{j+1} Y)]\otimes_{i+1} [(Z'\otimes_{j+1} W')\otimes_{1} (Z\otimes_{j+1} W)]
  \ar[dd]|{\eta^{1,j+1}_{X'Y'XY}\otimes_{i+1} \eta^{1,j+1}_{Z'W'ZW}}\\\\
  [(X'\otimes_{i+1} Z')\otimes_{1} (X\otimes_{i+1} Z)]\otimes_{j+1} [(Y'\otimes_{i+1} W')\otimes_{1} (Y\otimes_{i+1} W)]
  \ar[dd]|{\eta^{1,i+1}_{X'Z'XZ}\otimes_{j+1} \eta^{1,i+1}_{Y'W'YW}}
  &&[(X'\otimes_{1} X)\otimes_{j+1} (Y'\otimes_{1} Y)]\otimes_{i+1} [(Z'\otimes_{1} Z)\otimes_{j+1} (W'\otimes_{1} W)]
  \ar[dd]|{[M_{AA'A''}\otimes_{j+1} M_{BB'B''}]\otimes_{i+1} [M_{CC'C''}\otimes_{j+1} M_{DD'D''}]}
  \ar[ddll]|{\eta^{i+1,j+1}_{(X'\otimes_{1} X)(Y'\otimes_{1} Y)(Z'\otimes_{1} Z)(W'\otimes_{1} W)}}
  \\\\
  [(X'\otimes_{1} X)\otimes_{i+1} (Z'\otimes_{1} Z)]\otimes_{j+1} [(Y'\otimes_{1} Y)\otimes_{i+1} (W'\otimes_{1} W)]
  \ar[ddr]|{[M_{AA'A''}\otimes_{i+1} M_{CC'C''}]\otimes_{j+1} [M_{BB'B''}\otimes_{i+1} M_{DD'D''}]}
  &&[{\cal A}(A,A'')\otimes_{j+1} {\cal B}(B,B'')]\otimes_{i+1} [{\cal C}(C,C'')\otimes_{j+1} {\cal D}(D,D'')]
  \ar[ddl]|{\eta^{i+1,j+1}_{{\cal A}(A,A''){\cal B}(B,B''){\cal C}(C,C''){\cal D}(D,D'')}}\\\\
  &[{\cal A}(A,A'')\otimes_{i+1} {\cal C}(C,C'')]\otimes_{j+1} [{\cal B}(B,B'')\otimes_{i+1} {\cal D}(D,D'')]
  }
  $$
  \end{sideways}
  }
  \end{center}
  
  The lower quadrilateral commutes by naturality of $\eta$ and the upper hexagon commutes since it is an
  instance of the giant hexagonal interchange.
  \clearpage
  Again, as in the case of the same question for $\alpha^{(1)}$, the compatibility with the unit of $\eta^{(1)i,j}$ follows directly from the
  naturality of $\eta^{i+1,j+1}$ and the fact that $j_{[(A,B),(C,D)]} = [(j_A\otimes_{j+1} j_B)\otimes_{i+1} (j_C\otimes_{j+1} j_D)]$.
  
  Also as in the case of $\alpha^{(1)}$, the 2--naturality of $\eta^{(1)i,j}$ follows directly from the
  naturality of $\eta^{i+1,j+1}$ and the Mac Lane coherence theorem.
  
  Since $\alpha^{(1)}$ and $\eta^{(1)}$ are both defined in this way -- based upon $\alpha$ and $\eta$ -- we have 
  immediately that their
    ${\cal V}$--functor components satisfy all the axioms of the definition of a $k$--fold monoidal category. 
    At this level of course it is actually a $k$--fold
  monoidal 2--category.
  
  Notice that we have used all the axioms of a $k$--fold monoidal category. The external and internal unit conditions
  imply the unital nature of ${\cal V}$--Cat and the unit axioms for a product of ${\cal V}$--categories respectively.
  The external and internal associativities give us respectively the ${\cal V}$--functoriality of $\alpha^{(1)}$ and 
  the associativity of the composition morphisms for products of ${\cal V}$--categories. This reflects the dual nature
  of the latter two axioms that was pointed out in the braided case. Finally the giant hexagon gives us precisely the 
  ${\cal V}$--functoriality of $\eta^{(1)}.$ Notice also that we have used in each case the instance of the axiom 
  corresponding to $i=1; j=2..k.$  The remaining instances will be used as we iterate the categorical delooping. 
   \end{proof}

  \clearpage
  \newpage
   \MySection{Categories Enriched over ${\cal V}$--Cat}
    
  Here we generalize the definitions found in the appendix of \cite{Lyub} for a 
  ${\cal V}$--2--category. The main 
  difference is that we are considering ${\cal V}$ a $k$--fold monoidal category rather than symmetric monoidal. 
  I also choose to unpack the definition in terms of ${\cal V}$--functors, which the author of  \cite{Lyub}
  leaves implicit. Recall that the unit ${\cal V}$--category ${\cal I}$ has only one object $0$ and ${\cal I}(0,0)=I$ the 
  unit in ${\cal V}$. 
  
  \begin{definition} A (small,strict) ${\cal V}$--{\it 2--category} ${\bcal U}$ consists of
  
  \begin{enumerate}
      \item A set of objects $\left|{\bcal U}\right|$
      \item For each pair of objects $A,B \in \left|{\bcal U}\right|$ a ${\cal V}$--category ${\bcal U}(A,B).$
  
  Of course then ${\bcal U}(A,B)$ consists of a set of objects (which play the role of the 1--cells in a 2--category) 
  and for each pair $f,g \in \left|{\bcal U}(A,B)\right|$ an object 
  ${\bcal U}(A,B)(f,g) \in {\cal V}$ (which plays the role 
  of the hom--set of 2--cells in a 2--category.) Thus the vertical composition morphisms of these $hom_{2}$--objects are in ${\cal V}$:
  $$M_{fgh}:{\bcal U}(A,B)(g,h) \otimes_{1} {\bcal U}(A,B)(f,g) \to 
  {\bcal U}(A,B)(f,h)$$
  
  Also, the vertical identity for a 1-cell object $a \in \left|{\bcal U}(A,B)\right|$ is  $j_{a} : I \to {\bcal U}(A,B)(a,a)$.
  The associativity and the units of vertical composition are then those given by the respective axioms of enriched categories.  
      \item For each triple of objects $A,B,C \in 
  \left|{\bcal U}\right|$ a ${\cal V}$--functor
  $${\cal M}_{ABC}:{\bcal U}(B,C) \otimes^{(1)}_{1} {\bcal U}(A,B) \to {\bcal U}(A,C)$$ 
  Often I repress the subscripts. We denote ${\cal M}(h,f)$ as $hf$. 
  
  The family of morphisms indexed by pairs of objects $(g,f),(g',f') \in 
  \left|{\bcal U}(B,C) \otimes^{(1)}_{1} {\bcal U}(A,B)\right|$ furnishes the direct analogue of horizontal composition of 2-cells
  as can be seen by observing their domain and range in ${\cal V}$:
  $${\cal M}_{ABC_{(g,f)(g',f')}}:[{\bcal U}(B,C) \otimes^{(1)}_{1} 
  {\bcal U}(A,B)]((g,f),(g',f')) \to {\bcal U}(A,C)(gf,g'f')$$
  Recall that 
  $$[{\bcal U}(B,C) \otimes^{(1)}_{1} {\bcal U}(A,B)]((g,f),(g',f')) = {\bcal U}(B,C)(g,g') \otimes_{2} {\bcal U}(A,B)(f,f').$$
  We can now form the partial functors ${\cal M}(h,-):{\bcal U}(A,B) \to {\bcal U}(A,C)$ given by
  $$                             \xymatrix{
                                  {\bcal U}(A,B) = {\cal I} \otimes^{(1)}_{1} {\bcal U}(A,B)
                                  \ar[d]_{h \otimes^{(1)}_{1} 1}
                                  \\{\bcal U}(B,C) \otimes^{(1)}_{1} 
  {\bcal U}(A,B)
                                  \ar[d]_{\cal M}
                                  \\{\bcal U}(A,C)}.$$
  Where $h$ is here seen as the constant functor.
  
  Then $ {\cal M}(h,-)_{ff'}$  is given by
  $$\xymatrix{
                                  {\bcal U}(A,B)(f,f') = I \otimes_{2}
  {\bcal U}(A,B)(f,f')
                                  \ar[d]_{j_{h} \otimes_{2} 1}
                                  \\{\bcal U}(B,C)(h,h) \otimes_{2} {\bcal U}(A,B)(f,f')
                                  \ar[d]_{{\cal M}_{(h,f)(h,f')}}
                                  \\{\bcal U}(A,C)(hf,hf')}.$$
  This is the analogue of whiskering on the right. We can heuristically represent the objects of {\bcal U}(A,B) as arrows in 
  a diagram. The diagram for $ {\cal M}(h,-)_{ff'}$ should be
  $$
  
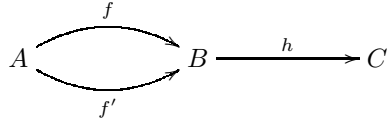
\xymatrix@R-=3pt{
  &\\
  A
  \ar@/^1pc/[rr]^f
  \ar@/_1pc/[rr]_{f'}
  &&B
  \ar[rr]^h
  &&C
  \\
  &\\
  }
  $$

  The other partial functors are ${\cal M}(-,f):{\bcal U}(B,C) \to {\bcal U}(A,C)$ given by
  $$                             \xymatrix{
                                  {\bcal U}(B,C) =  {\bcal U}(B,C) \otimes^{(1)}_{1} {\cal I}
                                  \ar[d]_{1 \otimes^{(1)}_{1} f}
                                  \\{\bcal U}(B,C) \otimes^{(1)}_{1} 
  {\bcal U}(A,B)
                                  \ar[d]_{\cal M}
                                  \\{\bcal U}(A,C)}.$$
  
  Then $ {\cal M}(-,f)_{hh'}$ is given by
  $$\xymatrix{
                                  {\bcal U}(B,C)(h,h') = 
  {\bcal U}(B,C)(h,h') \otimes_{2} I 
                                  \ar[d]_{1 \otimes_{2} j_{f}}
                                  \\{\bcal U}(B,C)(h,h') \otimes_{2} 
  {\bcal U}(A,B)(f,f)
                                  \ar[d]_{{\cal M}_{(h,f)(h',f)}}
                                  \\{\bcal U}(A,C)(hf,h'f)}.$$
  This is the analogue of whiskering on the left, as in 
  $$
  
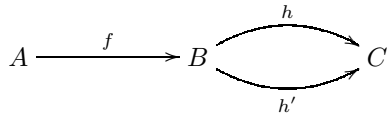
\xymatrix@R-=3pt{
  &&&\\
  A
  \ar[rr]^f
  &&B
  \ar@/^1pc/[rr]^h
  \ar@/_1pc/[rr]_{h'}
  &&C
  \\
  &&&\\
  }
  $$
  Notice that given any pair of partial functors, they cannot generally be combined to give a unique full functor since
  ${\cal V}$ is not symmetric.
     \item For each object $A \in \left|{\bcal U}\right|$ a ${\cal V}$--functor
  $${\cal J}_A: {\cal I} \to {\bcal U}(A,A)$$ 
  We denote ${\cal J}_A(0)$ as $1_{A}$. 
     \item (Associativity and unit axioms of a strict ${\cal V}$--2--category.) We require the pentagon and triangles of Definition
  \ref{V:Cat} to commute here as well. Since now the morphisms are 
  ${\cal V}$--functors this amounts to saying that the 
  functors given by the two legs of a diagram are equal. 
  For objects here we then have the equalities $(fg)h = f(gh)$ and $f1_{A} = f = 1_{B}f$
  
  For the  hom--object morphisms we have the following family of commuting diagrams for associativity, where the first bullet represents
  $$[({\bcal U}(C,D)\otimes^{(1)}_{1} {\bcal U}(B,C)) \otimes^{(1)}_{1} {\bcal U}(A,B)](((f,g),h),((f',g'),h'))$$
  and the reader may fill in the others
  $$
   
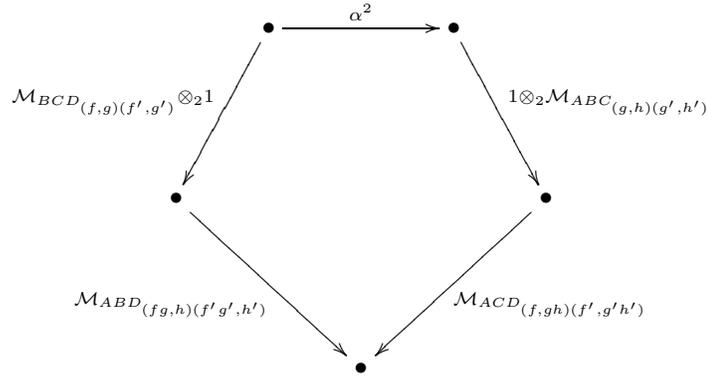
\xymatrix{
    &\bullet
    \ar[rr]^{\alpha^{2}}
    \ar[ddl]_{{\cal M}_{BCD_{(f,g)(f',g')}} \otimes_{2} 1}
    &&\bullet
    \ar[ddr]^{1 \otimes_{2} {\cal M}_{ABC_{(g,h)(g',h')}}}&\\\\
    \bullet
    \ar[ddrr]_{{\cal M}_{ABD_{(fg,h)(f'g',h')}}}
    &&&&\bullet
      \ar[ddll]^{{\cal M}_{ACD_{(f,gh)(f',g'h')}}}
      \\\\&&\bullet&&&
      }$$
      
   The heuristic diagram for this commutativity is
   $$
   
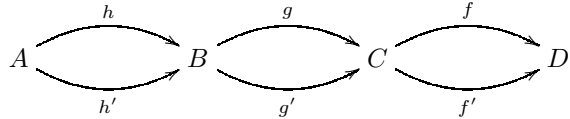
\xymatrix@R-=3pt{
   &&&\\
   A
   \ar@/^1pc/[rr]^h
   \ar@/_1pc/[rr]_{h'}
   &&B
   \ar@/^1pc/[rr]^g
   \ar@/_1pc/[rr]_{g'}
   &&C
   \ar@/^1pc/[rr]^f
   \ar@/_1pc/[rr]_{f'}
   &&D
   \\
   &&&\\
   }
  $$
      
   Some special cases in this family of commuting diagrams mentioned in \cite{Lyub} 
   are those described by the following heuristic diagrams.
   $$
    
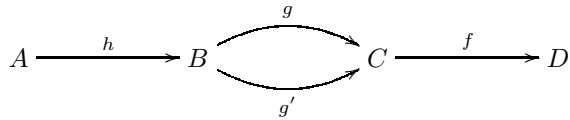
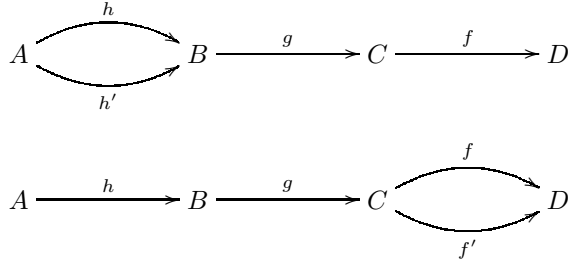
\xymatrix@R-=3pt{
    &&&\\
    A
    \ar[rr]^h
    &&B
    \ar@/^1pc/[rr]^g
    \ar@/_1pc/[rr]_{g'}
    &&C
    \ar[rr]^f
    &&D
    \\
    &&&\\
    }
  $$
  $$
   \xymatrix@R-=3pt{
   &&&\\
   A
   \ar@/^1pc/[rr]^h
   \ar@/_1pc/[rr]_{h'}
   &&B
   \ar[rr]^g
   &&C
   \ar[rr]^f
   &&D
   \\
   &&&\\
   }
  $$
  $$
   \xymatrix@R-=3pt{
   &&&\\
   A
   \ar[rr]^h
   &&B
   \ar[rr]^g
   &&C
   \ar@/^1pc/[rr]^f
   \ar@/_1pc/[rr]_{f'}
   &&D
   \\
   &&&\\
   }
  $$
   
   For the  unit morphisms we have that the triangles in the following diagram commute.

  $$
      \xymatrix@C=-5pt{
      [{\cal I}\otimes^{(1)}_1 {\bcal U}(A,B)]((0,f),(0,g))
      \ar[rrd]^{=}
      \ar[dd]_{{\cal J}_{B_{00}}\otimes_2 1}
      &&&&[{\bcal U}(A,B)\otimes^{(1)}_1 {\cal I}]((f,0),(g,0))
      \ar[dd]^{{1}\otimes_2 {\cal J}_{A_{00}}}
      \ar[lld]^{=}\\
      &&{\bcal U}(A,B)(f,g)\\
      [{\bcal U}(B,B)\otimes^{(1)}_1 {\bcal U}(A,B)]((1_B,f),(1_B,g))
      \ar[rru]|{{\cal M}_{ABB_{(1_B,f)(1_B,g)}}}
      &&&&[{\bcal U}(A,B)\otimes^{(1)}_1 {\bcal U}(A,A)]((f,1_A),(g,1_A))
      \ar[llu]|{{\cal M}_{AAB_{(f,1_A)(g,1_A)}}}
      }
   $$ 
   The heuristic diagrams for this commutativity are
   
   \noindent
                 		          \begin{center}
	          \resizebox{6in}{!}{
    $$
    \xymatrix@R-=3pt{
    &\ar@{=>}[dd]^{1_{1_A}}&&\\
    A
    \ar@/^1pc/[rr]^{1_A}
    \ar@/_1pc/[rr]_{1_A}
    &&A
    \ar@/^1pc/[rr]^f
    \ar@/_1pc/[rr]_g
    &&B
    \\
    &&&&\text{ }&\\
    }
    =
    \xymatrix@R-=3pt{
         &&\\
        A
        \ar@/^1pc/[rr]^f
        \ar@/_1pc/[rr]_g
        &&B
       \\
      &&&
      }
    =
    \xymatrix@R-=3pt{
      &&&\ar@{=>}[dd]^{1_{1_B}}\\
      A
      \ar@/^1pc/[rr]^f
      \ar@/_1pc/[rr]_g
      &&B
      \ar@/^1pc/[rr]^{1_B}
      \ar@/_1pc/[rr]_{1_B}
      &&B
      \\
      &&&\\
    }
  $$
  }
  \end{center}
  
   \end{enumerate}
    \end{definition}
    
   ${\cal V}$--functoriality of ${\cal M}$ and ${\cal J}$:
  First the ${\cal V}$--functoriality of ${\cal M}$ implies that the following (expanded) diagram commutes
  
  \noindent
              		          \begin{center}
	          \resizebox{6.5in}{!}{
  $$
   \xymatrix@C=-130pt{
   &({\bcal U}(B,C)(k,m)\otimes_1 {\bcal U}(B,C)(h,k))\otimes_2 ({\bcal U}(A,B)(g,l)\otimes_1 {\bcal U}(A,B)(f,g))
  \ar[rdd]^{M_{hkm}\otimes_2 M_{fgl}}
  \\\\
   ({\bcal U}(B,C)(k,m)\otimes_2 {\bcal U}(A,B)(g,l))\otimes_1 ({\bcal U}(B,C)(h,k)\otimes_2 {\bcal U}(A,B)(f,g))
  \ar[dd]^{{\cal M}_{ABC_{(k,g)(m,l)}}\otimes_1 {\cal M}_{ABC_{(h,f)(k,g)}}}
  \ar[ruu]^{\eta^{1,2}}
  &&{\bcal U}(B,C)(h,m)\otimes_2 {\bcal U}(A,B)(f,l)
  \ar[dd]^{{\cal M}_{ABC_{(h,f)(m,l)}}}
  \\\\
  {\bcal U}(A,C)(kg,ml)\otimes_1 {\bcal U}(A,C)(hf,kg)
  \ar[rr]^{M_{(hf)(kg)(ml)}}
  &&{\bcal U}(A,C)(hf,ml)
  \\
  }
  $$
  }
  \end{center}
  
  The heuristic diagram is
  $$
  
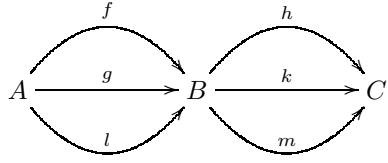
\xymatrix@R-=16pt{
  &&&
  \\
  A
  \ar@/^2pc/[rr]^f
  \ar[rr]^g
  \ar@/_2pc/[rr]^l
  && B
  \ar@/^2pc/[rr]^h
  \ar[rr]^k
  \ar@/_2pc/[rr]^m
  && C\\
  &&&&&
  }
  $$

  ${\cal V}$--functoriality of ${\cal M}$ implies ${\cal V}$--functoriality of the partial functors
  ${\cal M}(h,-).$ Special cases mentioned in \cite{Lyub} include those described by the diagrams
  $$
  
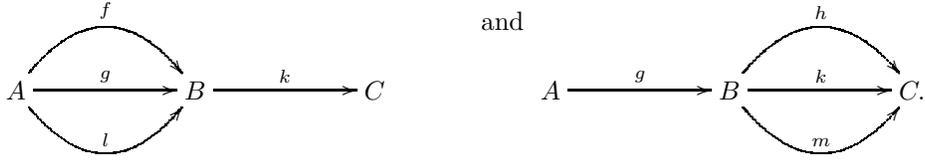
\xymatrix@R-=16pt{
  &&&
  \\
  A
  \ar@/^2pc/[rr]^f
  \ar[rr]^g
  \ar@/_2pc/[rr]^l
  && B
  \ar[rr]^k
  && C\\
  &&&&&
  }
  \text{ and }
  \xymatrix@R-=16pt{
  &&&
  \\
  A
  \ar[rr]^g
  && B
  \ar@/^2pc/[rr]^h
  \ar[rr]^k
  \ar@/_2pc/[rr]^m
  && C.\\
  &&&&&
  }
  $$
  Secondly the ${\cal V}$--functoriality of ${\cal M}$ implies that the following (expanded) diagram commutes
  $$
    
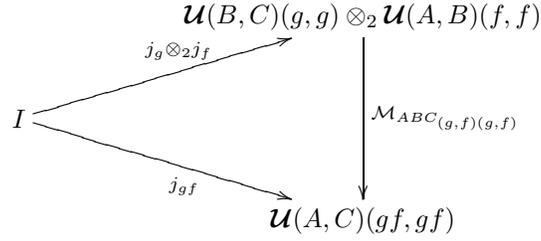
\xymatrix{
    &&{\bcal U}(B,C)(g,g)\otimes_2 {\bcal U}(A,B)(f,f)
    \ar[dd]^{{\cal M}_{ABC_{(g,f)(g,f)}}}\\
    I
    \ar[rru]^{j_{g}\otimes_2 j_{f}}
    \ar[rrd]_{j_{gf}}\\
    &&{\bcal U}(A,C)(gf,gf)
    }
   $$
  The heuristic diagram here is 
  $$
  
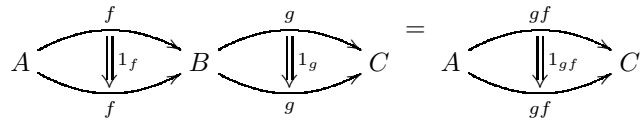
\xymatrix@R-=3pt{
  &\ar@{=>}[dd]^{1_f}&&\ar@{=>}[dd]^{1_g}\\
  A
  \ar@/^1pc/[rr]^f
  \ar@/_1pc/[rr]_f
  &&B
  \ar@/^1pc/[rr]^g
  \ar@/_1pc/[rr]_g
  &&C
  \\
  &&&\\
  }
  =
  \xymatrix@R-=3pt{
  &\ar@{=>}[dd]^{1_{gf}}&&\\
  A
  \ar@/^1pc/[rr]^{gf}
  \ar@/_1pc/[rr]_{gf}
  &&C
  \\
  &&&\\
  } 
  $$
  In addition, the ${\cal V}$--functoriality of ${\cal J}$ implies that the following (expanded) diagram commutes
  $$
    
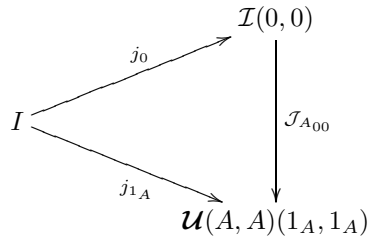
\xymatrix{
    &&{\cal I}(0,0)
    \ar[dd]^{{\cal J}_{A_{00}}}\\
    I
    \ar[rru]^{j_{0}}
    \ar[rrd]_{j_{1_A}}\\
    &&{\bcal U}(A,A)(1_A,1_A)
    }
  $$
  Which means that 
  $${\cal J}_{A_{00}}: I \to {\bcal U}(A,A)(1_{A},1_{A}) = j_{1_{A}}.$$
  In other words the ``horizontal'' unit for the object $1_A$ is the same as the ``vertical'' unit for $1_A.$

  I now describe the (strict) 3--category  ${\cal V}$--2--Cat (or ${\cal V}$--Cat--Cat)
  whose objects are (strict, small)  ${\cal V}$--2--categories.
    We are guided by the definitions of  ${\cal V}$--functor and  ${\cal V}$--natural transformation 
    as well as by the definitions of 2--functor, 2--natural transformation, and modification.
    
    \begin{definition} \label{2:enriched:funct}
    For two ${\cal V}$--2--categories ${\bcal U}$ and ${\bcal W}$ 
    a ${\cal V}${\it --2--functor} $T:{\bcal U} \to {\bcal W}$ is a function on objects 
    $\left|{\bcal U}\right| \to \left|{\bcal W}\right|$ and a family of ${\cal V}$--functors
    $T_{UU'}:{\bcal U}(U,U') \to {\bcal W}(TU,TU').$ These latter obey commutativity of the usual diagrams.
    \begin{enumerate}
     \item For $U,U',U'' \in \left|{\bcal U}\right|$
     $$
       \xymatrix@C=65pt@R=65pt{
       &\bullet
       \ar[rr]^{{\cal M}_{UU'U''}}
       \ar[d]^{T_{U'U''} \otimes^{(1)}_1 T_{UU'}}
       &&\bullet
       \ar[d]^{T_{UU''}}&\\
       &\bullet
       \ar[rr]_{{\cal M}_{(TU)(TU')(TU'')}}
       &&\bullet
       }
     $$
     \item
     $$
       \xymatrix{
       &&\bullet
       \ar[dd]^{T_{UU}}\\
       {\cal I}
       \ar[rru]^{{\cal J}_U}
       \ar[rrd]_{{\cal J}_{TU}}\\
       &&\bullet
       }
     $$
     
     \end{enumerate}
    \end{definition}
     For objects this means that $T_{U'U''}(f)T_{UU'}(g) = T_{UU''}(fg)$ and $T_{UU}(1_U) = 1_{TU}.$
     The reader should unpack both diagrams into terms of hom--object morphisms and ${\cal V}$--functoriality.
     Composition of ${\cal V}$--2--functors is just composition of functions and components.

    
    \begin{definition} \label{2:enr:nat:trans}
    A ${\cal V}${\it --2--natural transformation} $\alpha:T\to S:{\bcal U}\to {\bcal W}$ is a function sending
    each $U \in  \left|{\bcal U}\right|$ to a 
    ${\cal V}$--functor $\alpha_{U}: {\cal I}\to {\bcal W}(TU,SU)$ in such a way that we have commutativity of
    
    $$
      \xymatrix{
      &{\cal I} \otimes^{(1)}_1 {\bcal U}(U,U')
      \ar[rr]^-{ \alpha_{U'} \otimes^{(1)}_1 T_{UU'}}
      &&{\bcal W}(TU',SU') \otimes^{(1)}_1 {\bcal W}(TU,TU')
      \ar[rd]^-{{\cal M}}
    \\
      {\bcal U}(U,U')
      \ar[ru]^{=}
    \ar[rd]_{=}
      &&&&{\bcal W}(TU,SU')
    \\
      &{\bcal U}(U,U') \otimes^{(1)}_1 {\cal I}
      \ar[rr]_-{S_{UU'} \otimes^{(1)}_1 \alpha_{U}}
      &&{\bcal W}(SU,SU') \otimes^{(1)}_1 {\bcal W}(TU,SU)
      \ar[ru]^-{{\cal M}}
      }
    $$
    
    \end{definition}
   Unpacking this a bit, we see that $\alpha_U$ is an object $q = \alpha_U(0)$ in 
   the ${\cal V}$--category ${\bcal W}(TU,SU)$ and 
   a morphism $\alpha_{U_{00}}: I \to {\bcal W}(TU,SU)(q,q).$ By the ${\cal V}$--functoriality of $\alpha_U$
   we see that $\alpha_{U_{00}} = j_q.$ The axiom then states that $q'T_{UU'}(f) = S_{UU'}(f)q$ for all $f$,
   and that 
   $${\cal M}_{(TU)(TU')(SU')_{(q',T_{UU'}(f))(q',T_{UU'}(g))}} \circ (j_{q'} \otimes_2 T_{UU'_{fg}})
   ={\cal M}_{(TU)(SU)(SU')_{(S_{UU'}(f),q)(S_{UU'}(g),q)}} \circ (S_{UU'_{fg}}\otimes_2 j_q )$$
   This is directly analogous to the usual definition of 2--natural transformation by whisker diagrams.
   
   Vertical composition of ${\cal V}$--2--natural transformations is as expected.
   $(\beta \circ \alpha)_U =$
   $$
   \xymatrix{
   {\cal I} \otimes^{(1)}_1 {\cal I}
   \ar[d]^{\beta_U \otimes^{(1)}_1 \alpha_U}
   \\
   {\bcal W}(SU,RU) \otimes^{(1)}_1 {\bcal W}(TU,SU)
   \ar[d]^{{\cal M}}
   \\
   {\bcal W}(TU,RU)
   }
   $$
   Identity 2--cells for vertical composition are ${\cal V}$--2--natural transformations $\mathbf{1}_T:T\to T$ where 
   $(\mathbf{1}_T)_U = {\cal J}_{TU}.$
   Left and right whiskering of ${\cal V}$--2--functors onto ${\cal V}$--2--natural transformations are given by precisely the same 
   descriptions as in the low dimensional case, with $I$ replaced by ${\cal I}$, etc.

      \begin{definition} \label{modification}
      Given two ${\cal V}$--2--natural transformations a 
    ${\cal V}${\it --modification} between them $\mu:\theta\to \phi :T \to S:{\bcal U}\to {\bcal W}$
    is a function that sends each object $U \in \left|{\bcal U}\right|$ to a morphism
    $\mu_U:I \to {\bcal W}(TU,SU)(\theta_U(0),\phi_U(0))$ in such a way that the following diagram commutes.
    (Let $\theta_U(0) = q$, $\phi_U(0) = \hat{q}$, $\theta_{U'}(0) = q'$ and $\phi_{U'}(0) = \hat{q'}.$)
    
     \noindent
                 		          \begin{center}
	          \resizebox{6.5in}{!}{
     $$
        \xymatrix@C=-35pt@R=20pt{
        &
        &{\bcal W}(TU',SU')(q',\hat{q'}) \otimes_2 {\bcal W}(TU,TU')(T_{UU'}(f),T_{UU'}(g))
        \ar[rd]^-{{\cal M}}
      \\
      &
      I \otimes_2 {\bcal U}(U,U')(f,g)
        \ar[ru]_{\text{ }\mu_{U'} \otimes_2 T_{UU'_{fg}}}
      &&{\bcal W}(TU,SU')(q'T_{UU'}(f),\hat{q'}T_{UU'}(g))
      \ar@{=}[dd]
      \\
      {\bcal U}(U,U')(f,g)
        \ar[ru]^{=}
      \ar[rd]_{=}
      \\
      &
      {\bcal U}(U,U')(f,g) \otimes_2 I
        \ar[rd]^{\text{ }S_{UU'_{fg}} \otimes_2 \mu_{U}}
      &&{\bcal W}(TU,SU')(S_{UU'}(f)q,S_{UU'}(g)\hat{q})
      \\
        &
        &{\bcal W}(SU,SU')(S_{UU'}(f),S_{UU'}(g)) \otimes_2 {\bcal W}(TU,SU)(q,\hat{q})
        \ar[ru]^-{{\cal M}}
        }
      $$
 }
 \end{center}
 
 \end{definition}

      This is directly analogous to the usual definition of modification described in section 1.
    Notice that since $\theta_{U_{00}} = j_{\theta_U(0)}$ for all ${\cal V}$--2--natural transformations $\theta$ we have 
    that the morphism $\mu_U$ seen as a ``family'' consisting of a single morphism 
    (corresponding to $0 \in \left|{\cal I}\right|$)
    constitutes a ${\cal V}$--natural transformation from $\theta_U$ to $\phi_U.$ ``Vertical'' compositions of 
    modifications are given by
    the compositions of these underlying ${\cal V}$--natural transformations as described in section 1. 
    Thus identities  $\mathbf{1}_{\alpha}$ for this
    composition are families of ${\cal V}$--natural equivalences. Since $\alpha_U$ is a 
    ${\cal V}$--functor from ${\cal I}$ to ${\bcal W}(TU,SU)$ this means specifically that $((\mathbf{1}_{\alpha})_U)_0 = j_{\alpha_{U}(0)} = j_q.$
      
\begin{theorem}
${\cal V}$--2--categories, ${\cal V}$--2--functors, ${\cal V}$--2--natural transformations
and ${\cal V}$--modifications form a 3--category called ${\cal V}$--2--Cat. 
\end{theorem}
The proof for this is long and
tedious or long and interesting, depending on your point of view. Since I am in the latter camp, the details are spelled out in
the forthcoming \cite{Forcey}. That paper also includes the definitions of ${\cal V}$--$n$--categories and of the
morphisms of ${\cal V}$--$n$--Cat. 

    For ${\cal V}$ $k$--fold monoidal we have demonstrated that ${\cal V}$--Cat is $(k-1)$--fold monoidal. It is straightforward to show
    that this process continues, i.e. that ${\cal V}$--2--Cat is $(k-2)$--fold monoidal. We need a unit ${\cal V}$--2--category; this is easily 
  given as ${\bcal I}$ where $\left|{\bcal I}\right| = \{\mathbf{0}\}$ and ${\bcal I}(\mathbf{0},\mathbf{0}) = {\cal I}.$
  Products of ${\cal V}$--2--categories are given by ${\bcal U} \otimes^{(2)}_i {\bcal W}$ for $i=1...k-2.$ Objects are pairs of objects as usual,
  and that there are exactly $k-2$ products is seen when the definition of hom--objects is given. In ${\cal V}$--2--Cat,
  $$ [{\bcal U} \otimes^{(2)}_i {\bcal W}]((U,W),(U',W')) = {\bcal U}(U,U') \otimes^{(1)}_{i+1} {\bcal W}(W,W') $$
  Thus we have that
  $$ [{\bcal U} \otimes^{(2)}_i {\bcal W}]((U,W),(U',W'))((f,f'),(g,g')) $$
  $$ = [{\bcal U}(U,U') \otimes^{(1)}_{i+1} {\bcal W}(W,W')]((f,f'),(g,g')) $$
  $$ = {\bcal U}(U,U')(f,g) \otimes_{i+2} {\bcal W}(W,W')(f',g') $$
  Of course we now consider enrichment over ${\cal V}$--Cat with composition given by ${\cal M}$
  and units given by ${\cal J}.$
  
  The definitions of $\alpha^{(2)i}$ and $\eta^{(2)i,j}$  are just as in the lower dimensional
  case. 
  For instance, $\alpha^{(2)i}$ will now be a 3--natural transformation, that is, a family of
  ${\cal V}$--2--functors 
  $$\alpha^{(2)i}_{{\bcal U}{\bcal V}{\bcal W}}:({\bcal U}\otimes^{(2)}_i {\bcal V}) \otimes^{(2)}_i {\bcal W} \to {\bcal U}\otimes^{(2)}_i ({\bcal V}\otimes^{(2)}_i {\bcal W}).$$
  To each of these is associated a family of ${\cal V}$--functors 
  $$\alpha^{(2)i}_{{\bcal U}{\bcal V}{\bcal W}_{(U,V,W)(U',V',W')}} = \alpha^{(1)i+1}_{{\bcal U}(U,U'){\bcal V}(V,V'){\bcal W}(W,W')}$$
  to each of which is associated a family of hom--object morphisms
  $$\alpha^{(2)i}_{{\bcal U}{\bcal V}{\bcal W}_{(U,V,W)(U',V',W')_{(f,g,h)(f',g',h')}}} = \alpha^{i+2}_{{\bcal U}(U,U')(f,f'){\bcal V}(V,V')(g,g'){\bcal W}(W,W')(h,h')}.$$
  Verifications that these define
  a valid $(k-2)$--fold monoidal 3--category all follow just as in the lower dimensional
  case. The facts about the hom--object morphisms are shown by previous demonstrations exactly.
  The facts about the ${\cal V}$--functors are shown in an analogous way,
  but now using the original $k$--fold monoidal category axioms that involve i=2.

}
\end{document}